%% file: main.tex
\documentclass[final,hidelinks,onefignum,onetabnum]{siamart220329}

\usepackage{cancel}
\input{ex_shared}

\usepackage{bm}
\usepackage{subcaption}

\usepackage{mathtools}
\usepackage{lineno}

\usepackage{calligra,amsmath,amssymb}

\begin{document}

\maketitle

\begin{abstract}
This work introduces a parallel and rank-adaptive matrix integrator for dynamical low-rank approximation. The method is related to the previously proposed rank-adaptive basis update \& Galerkin (BUG) integrator but differs significantly in that all arising differential equations, both for the basis \emph{and} the Galerkin coefficients, are solved in parallel. Moreover, this approach eliminates the need for a potentially costly coefficient update with augmented basis matrices. The integrator also incorporates a new step rejection strategy that enhances the robustness of both the parallel integrator and the BUG integrator. By construction, the parallel integrator inherits the robust error bound of the BUG and projector-splitting integrators. Comparisons of the parallel and BUG integrators are presented by a series of numerical experiments which demonstrate the efficiency of the proposed method, for problems from radiative transfer and radiation therapy.
\end{abstract}

\begin{keywords}
Dynamical low-rank approximation, rank adaptivity, time integration
\end{keywords}

\begin{MSCcodes}
68Q25, 68R10, 68U05
\end{MSCcodes}

\section{Introduction}
Dynamical low-rank approximation (DLRA)~\cite{koch2007dynamical} has shown to be an efficient approach to reduce the computational cost and memory footprint of solving large matrix differential equations
\begin{equation}\label{ode}
\dot \bfA (t) = \bfF(t,\bfA(t)) \, .
\end{equation}
The core building block of dynamical low-rank approximation is to evolve an approximate solution on the low-rank manifold $\mathcal{M}$ by orthogonally projecting the vector-field $\bfF(t,\bfA(t))$ onto the tangent space of $\mathcal{M}$ at each instance in time. This approach is known as the time-dependent variational principle for the approximation manifold $\mathcal{M}$ in the physical and chemical literature, going back to Dirac \cite{dirac1930}  in the context of a quantum-mechanical problem that is approximated by separation of variables.

This approach evolves entirely on the low-rank manifold and thus differs from schemes which leave the manifold and then retract the numerical approximation back to the manifold.  However, a severe challenge to the use of numerical time marching schemes is posed by high curvature of the manifold along the solution. The curvature is proportional to the inverse of the smallest nonzero singular value of the numerical approximation and the exact solution~\cite[\S 4]{koch2007dynamical}. The smallest nonzero singular value is typically very small, because the accuracy of a low-rank approximation corresponds to the largest discarded singular value of the solution, and then the smallest retained singular value is usually not much larger. The smallest nonzero singular value is therefore the smaller the higher the rank and the higher the accuracy requirements.
%

Overcoming the challenge posed by small singular values in constructing robust numerical schemes for dynamical low-rank approximation has been a major difficulty. So far, two numerical integrators have been developed that resolve this issue: the projector-splitting integrator~\cite{LubO14} and the basis update \& Galerkin (BUG) integrator~\cite{ceruti2022unconventional}. The projector-splitting integrator applies a splitting scheme to 
the orthogonal projection onto the tangent space, which is an alternating sum of subprojections; it has been proved that this integrator is not affected by the presence of small singular values~\cite{kieri2016robust}. The BUG integrator, along with its rank-adaptive extension~\cite{ceruti2022rank}, updates the basis matrices similarly to the projector-splitting integrator but differs by performing the basis update in parallel while not requiring a backward-in-time step, which is important in strongly dissipative problems. The robustness of the BUG integrator is achieved through a perturbation argument applied to properties of the projector-splitting integrator.
The basis augmentation step proposed in~\cite{ceruti2022rank} allows for rank-adaptivity while exhibiting remarkable preservation properties. 

Here we address some issues that remain with the rank-adaptive BUG integrator: First, while the basis update can be computed in parallel, the Galerkin step must be computed in a subsequent serial step. Second, the Galerkin update of double dimension requires a potentially expensive projection of the vector field onto the subspace spanned by the old \emph{and} new bases. Third,
the error due to the tangent space projection, which depends on the chosen rank, is not controlled in~\cite{ceruti2022rank}.
%

Already in \cite[\S 3.1]{ceruti2022unconventional}, it has been noted that BUG integrators offer the potential to compute the basis update \emph{and} the Galerkin step in parallel. However, the modification given there lacks the robustness of the serial Galerkin step, as it requires a matrix inverse that is not under control except for step-sizes proportional to the smallest singular value. Another attempt to achieve full parallelism is found in the context of quantum physics \cite{dolgov2020parTVDP}, where the projector-splitting integrator is employed. The time-evolution of matrix product states is investigated and a computational algorithm using the message passing interface (MPI) architecture is achieved by recursively halving the matrix product state. Although this approach offers parallel features, it also suffers from small singular values~\cite[\S III.D]{dolgov2020parTVDP}. 
%
A different approach exploits parallelism in time as proposed in \cite{carrel2023parareal}, which combines parallel-in-time methods~\cite{gander2007parareal} with dynamical low-rank approximation techniques. 
%


In this work, we introduce a novel rank-adaptive integrator which updates the bases \emph{and} performs the Galerkin update in parallel. Unlike the rank-adaptive BUG integrator~\cite{ceruti2022rank}, it does not require a Galerkin update with the augmented basis, which reduces computational costs. The new parallel integrator still shares the robust error bound of the projector-splitting and BUG integrators while no longer fulfilling the exactness property. The construction of the integrator is based on identifying new and old basis information in the rank-adaptive BUG integrator while using the fact that neglecting serially computed contributions still preserves the asymptotic error of the original integrator. 
Moreover, a step rejection strategy is proposed that can be used to efficiently control the low-rank error. This strategy applies to the parallel as well as the rank-adaptive BUG integrator of \cite{ceruti2022rank}.

It should be noted that while the new parallel integrator is presented here for the matrix case of dynamical low-rank approximation, it is expected that its concepts extend to tensor trains and general tree tensor networks, where a much more substantial reduction of computation times by parallelism is to be gained. This will be explored in detail in future work, which will be based on ideas in the notationally and technically less demanding matrix case presented here.

The paper is structured as follows: In Section 2, we recapitulate the BUG integrator of \cite{ceruti2022unconventional} and its rank-adaptive formulation~\cite{ceruti2022rank}.
In Section 3, the new fully parallel integrator together with its rationale is introduced; a new step rejection strategy is discussed. In Section 4, the error analysis of the parallel integrator is presented. 
In Section 5, numerical experiments with examples from radiative transfer and radiation therapy validate the theoretical results of the new fully parallel integrator.

\section{Recap: the rank-adaptive BUG integrator of \cite{ceruti2022rank}}
\label{sec:rabug}
The rank-adaptive parallel integrator presented in this paper has conceptual similarities to the Basis Update \& Galerkin (BUG) integrator of \cite{ceruti2022unconventional} and its rank-adaptive version \cite{ceruti2022rank}. The rank-adaptive BUG integrator has been shown to be robust to small singular values in the exact solution $\bfA(t)$ or its low-rank approximation, and it has good conservation properties for Hamiltonian systems and dissipation properties for gradient systems up to the truncation tolerance $\vartheta$.
To understand how the new parallel integrator relates to and differs from these previously proposed integrators, it is useful to recall the rank-adaptive BUG integrator of \cite{ceruti2022rank}. Furthermore, in Section~\ref{subsec:rejection} we introduce criteria for step rejection, which are useful for the rank-adaptive integrator of \cite{ceruti2022rank} and play an even more important role for the parallel integrator of this paper.

\subsection{Formulation of the rank-adaptive BUG integrator}
	A time step of integration from time $t_0$ to $t_1=t_0+h$  starts from a factored rank-$r_0$ matrix 
	$\bfY_0=\bfU_0\bfS_0\bfV_0^\top$ and computes an updated factorization $\bfY_1=\bfU_1\bfS_1\bfV_1^\top$ of rank $r_1 \le 2r_0$.  We let $r=r_0$ and we put a hat on quantities related to the augmented rank $\widehat r$ with $r \le \wh r \le 2r$ (and typically $\wh r= 2r$). Let $\vartheta$ be a given rank truncation tolerance, just as a local error tolerance is provided by the user of an integrator with adaptive stepsizes.
	
	\begin{enumerate}
		\item 
		Compute augmented basis matrices $\wh \bfU\in \R^{m\times \widehat r}$ and $\wh \bfV\in \R^{n\times \widehat r}$ (in parallel):
		\\[2mm]
		\textbf{K-step}:
		Integrate from $t=t_0$ to $t_1$ the $m \times r$ matrix differential equation
		\begin{equation}\label{K-step-rabug} 
		\dot{\textbf{K}}(t) = \bfF(t, \textbf{K}(t) \bfV_0^\top) \bfV_0, \qquad \textbf{K}(t_0) = \bfU_0 \bfS_0.
		\end{equation}
		Determine the columns of $\wh \bfU\in \R^{m\times \widehat r}$ with $\wh r \le 2r$ as an orthonormal basis of the range 
		of the $m\times 2r$ matrix $(\textbf{K}(t_1),\bfU_0)$ (e.g.~by QR decomposition)
		and compute the $\widehat r\times r$ matrix $\wh \bfM= \wh\bfU^\top \bfU_0$.
		\\[2mm]
		\textbf{L-step}: 
		Integrate from $t=t_0$ to $t_1$ the $n \times r$ matrix differential equation
		\begin{equation}\label{L-step-rabug} 
		\dot{\textbf{L}}(t) =\bfF(t, \bfU_0 \textbf{L}(t)^\top)^\top  \bfU_0, \qquad \textbf{L}(t_0) = \bfV_0 {\bfS}_0^\top. 
		\end{equation}
		Determine the columns of $\wh \bfV\in \R^{n\times \wh r}$ as an orthonormal basis of the range 
		of the $n\times 2r$ matrix $(\textbf{L}(t_1),\bfV_0)$ (e.g.~by QR decomposition)
		and compute the $\wh r\times r$ matrix $\wh \bfN= \wh\bfV^\top \bfV_0$.
		\\[-2mm]
		\item
		\textbf{Augment and update} ${\bfS}_0 \rightarrow {\wh\bfS}(t_1)$\,: \\[1mm]
		\textbf{S-step}:  Integrate from $t=t_0$ to $t_1$ the $\wh r \times \wh r$ matrix differential equation
		\begin{equation}\label{S-step-rabug} 
		\dot{\wh\bfS}(t) =  \wh\bfU^\top \bfF(t, \wh\bfU \wh\bfS(t) \wh\bfV^\top) \wh\bfV, 
		\qquad \wh\bfS(t_0) = \wh \bfM \bfS_0 \wh \bfN^\top.
		\end{equation}
		\item \textbf{Truncate}:
		Compute the SVD $\; \wh\bfS(t_1)= \wh \bfP \wh \bfSigma \wh \bfQ^\top$ with $\wh\bfSigma=\diag(\sigma_j)$ and truncate to the tolerance~$\vartheta$: Choose the new rank $r_1\le \wh r \le 2r$ as the 
		minimal number $r_1$ such that 
		$$
		\biggl(\ \sum_{j=r_1+1}^{\wh r} \sigma_j^2 \biggr)^{1/2} \le \vartheta.
		$$
		Compute the new factors for the approximation of $\bfY(t_1)$ as follows:
Let $\bfS_1$ be the $r_1\times r_1$ diagonal matrix with the $r_1$ largest singular values and let $\bfP_1\in \R^{\wh r\times r_1}$ and $\bfQ_1\in \R^{\wh r\times r_1}$ contain the first $r_1$ columns of $\wh \bfP$ and $\wh \bfQ$, respectively. Finally, set $\bfU_1 = \wh \bfU \bfP_1\in \R^{m\times r_1}$ and
$\bfV_1 = \wh \bfV \bfQ_1 \in \R^{n\times r_1}$.
	\end{enumerate} 
	The approximation after one time step is given by 
	\begin{equation}\label{Y1-rabug}
	 \bfY_1 = \bfU_1 \bfS_1 \bfV_1^\top \approx \bfY(t_1).
	 \end{equation}
	Then,  $\bfY_1$ is taken as the starting value for the next step, which computes $\bfY_2$ in factorized form, etc. 
	
The $m\times r$, $n\times r$ and $\wh r\times \wh r$ matrix differential equations in the substeps are solved approximately using a standard integrator, such as a Runge--Kutta method or an exponential integrator when $\bfF$ is predominantly linear.

To simplify the notation, we ignored that there might be different augmented ranks $\wh r_U = \text{rank}(\textbf{K}(t_1),\bfU_0)$ and $\wh r_V = \text{rank}(\textbf{L}(t_1),\bfV_0)$. This is readily taken into account in the above algorithm. In any case, typically one has $\wh r=\wh r_U= \wh r_V = 2r$.

\subsection{Interpretation}
The K-step updates and augments the left basis from  $\bfU_0$ to $\wh \bfU$. The differential equation for $\bfK$ can be viewed as a Galerkin method in the subspace of matrices $\{ \bfK \bfV_0^\top\,:\, \bfK \in \R^{m\times r}\} \subset \R^{m\times n}$ with fixed corange $\text{Ran}(\bfV_0)$. 
Similarly, the differential equation for $\bfL$ is a Galerkin method in the subspace of matrices $\{ \bfU_0 \bfL^\top\,:\, \bfL \in \R^{n\times r}\} \subset \R^{m\times n}$ with fixed range $\text{Ran}(\bfU_0)$.

The S-step is a Galerkin method for the differential equation \eqref{ode} in the space of matrices $\wh\bfU \bfS \wh\bfV^\top$ generated by the extended basis matrices $\wh\bfU$ and $\wh\bfV$, with the projected starting value 
$\wh \bfY_0 =\wh \bfU \wh \bfU^\top \bfY_0 \wh \bfV \wh \bfV^\top = \wh\bfU \wh \bfS(t_0) \wh \bfV^\top$. 
Since $\text{Ran}\,\bfU_0\subset \text{Ran}\,\wh\bfU$ and $\text{Ran}\,\bfV_0\subset \text{Ran}\,\wh\bfV$, the starting value is not changed: $\wh \bfY_0= \bfY_0=\bfU_0\bfS_0\bfV_0^\top$.

\section{The fully parallel integrator}
In the following, we present a novel parallel integrator that shares the robust error bound of the projector-splitting integrator \cite{LubO14} and of the Basis Update \& Galerkin (BUG) integrator \cite{ceruti2022unconventional} as well as its rank-adaptive extension \cite{ceruti2022rank}. Important properties of the integrator are:
\begin{itemize}
    \item All arising differential equations (for $\bfK$, $\bfL$ and $\bfS$) can be solved in parallel.
    \item All differential equations are evolved forward in time. This property is also valid for the BUG integrator, however not for the projector-splitting integrator, which requires a backward-in-time S-step.
    \item The integrator is rank-adaptive. That is, ranks are chosen in an automated way based on a user-determined tolerance parameter. Compared to the rank-adaptive BUG integrator, the new integrator does not require the sequential update of the $2r\times 2r$ augmented coefficient matrix $\widehat S$, which by itself is more costly than the parallel $r\times r$ S-step given below.
\end{itemize}
In contrast to the projector-splitting and the BUG integrators, the parallel integrator does not have the exactness property. That is, when the exact solution $\bfA(t)$ is of rank $r$ and the right-hand side does not depend on the solution, i.e., $\bfF = \bfF(t)=\dot\bfA(t)$, the solution computed with the parallel integrator is in general not reproduced exactly. Moreover, unlike the rank-adaptive BUG integrator, the new integrator does not preserve the energy (up to a multiple of the truncation parameter) for Hamiltonian systems and does not necessarily diminish the energy (up to a multiple of the truncation parameter) for gradient systems.
\subsection{Formulation of the algorithm}
    A time step from $t_0$ to $t_1$ where the numerical solution at time $t_0$ is given by the factored rank-$r$ matrix $\bfY_0 = \bfU_0 \bfS_0 \bfV_0^{\top}$ reads:
	\begin{enumerate}
		\item 
		Determine the augmented basis matrices $\wh \bfU\in \R^{m\times 2r}$ and $\wh \bfV\in \R^{n\times 2r}$ as well as the coefficient matrix $\bar{\bfS}(t_1) \in\R^{r \times r}$ (in parallel):
		\\[2mm]
		\textbf{K-step}:
		Integrate from $t=t_0$ to $t_1$ the $m \times r$ matrix differential equation
		\begin{equation}\label{K-step} 
		\dot{\textbf{K}}(t) = \bfF(t, \textbf{K}(t) \bfV_0^\top) \bfV_0, \qquad \textbf{K}(t_0) = \bfU_0 \bfS_0.
		\end{equation}
		Determine $\wh \bfU = (\bfU_0, \widetilde \bfU_1) \in \R^{m\times 2r}$ as an orthonormal basis of the range 
		of the $m\times 2r$ matrix $(\bfU_0, \textbf{K}(t_1))$ (e.g.~by QR decomposition), where $ \widetilde \bfU_1 \in \R^{m\times r}$ 
		is filled with zero columns if $(\bfU_0, \textbf{K}(t_1))$ has rank less than $2r$.\\
		Compute the matrix $\mathbf{\widetilde S}_1^K = \widetilde\bfU_1^{\top}\bfK(t_1)\in\R^{r\times r}$.
		\\[2mm]
		\textbf{L-step}: 
		Integrate from $t=t_0$ to $t_1$ the $n \times r$ matrix differential equation
		\begin{equation}\label{L-step} 
		\dot{\textbf{L}}(t) =\bfF(t, \bfU_0 \textbf{L}(t)^\top)^\top  \bfU_0, \qquad \textbf{L}(t_0) = \bfV_0 {\bfS}_0^\top. 
		\end{equation}
		Determine $\wh \bfV = (\bfV_0, \widetilde \bfV_1) \in \R^{n\times 2r}$ as an orthonormal basis of the range 
		of the $n\times 2r$ matrix $(\bfV_0, \textbf{L}(t_1))$ (e.g.~by QR decomposition), where $ \widetilde \bfV_1 \in \R^{n\times r}$ 
		is filled with zero columns if $(\bfV_0, \textbf{L}(t_1))$ has rank less than $2r$.\\
		Compute the matrix $\mathbf{\widetilde S}_1^L = \bfL(t_1)^{\top}\widetilde\bfV_1\in\R^{r\times r}$.
		\\[2mm]
		\textbf{S-step}:
		Integrate from $t=t_0$ to $t_1$ the $r \times r$ matrix differential equation
		\begin{equation}\label{S-step} 
		\dot{\bar\bfS}(t) =  \bfU_0^\top \bfF(t, \bfU_0 \bar\bfS(t) \bfV_0^\top) \bfV_0, 
		\qquad \bar\bfS(t_0) = \bfS_0.
		\end{equation}
		\item \textbf{Augment}:
		Set up the augmented coefficient matrix $\widehat \bfS_1 \in\R^{2r\times 2r}$ as
		\begin{align} \label{wh-S1}
            \widehat \bfS_1 = 
            \begin{pmatrix} 
                \bar \bfS(t_1) & \widetilde \bfS_1^L \\
                \widetilde \bfS_1^K & \mathbf{0}
            \end{pmatrix}.
        \end{align}
        \item \textbf{Truncate}:
        Compute the singular value decomposition $\; \wh\bfS_1= \wh \bfP \wh \bfSigma \wh \bfQ^\top$ where $\wh\bfSigma=\diag(\sigma_j)$. Truncate to the tolerance~$\vartheta$ by choosing the new rank $r_1\le 2r$ as the 
		minimal number $r_1$ such that 
		\begin{equation}\label{vartheta}
		\biggl(\ \sum_{j=r_1+1}^{2r} \sigma_j^2 \biggr)^{1/2} \le \vartheta.
		\end{equation}
		Compute the new factors for the approximation of $\bfA(t_1)$ as follows:
Let $\bfS_1$ be the $r_1\times r_1$ diagonal matrix with the $r_1$ largest singular values and let $\bfP_1\in \R^{2r\times r_1}$ and $\bfQ_1\in \R^{2r\times r_1}$ contain the first $r_1$ columns of $\wh \bfP$ and $\wh \bfQ$, respectively. Finally, set $\bfU_1 = \wh \bfU \bfP_1\in \R^{m\times r_1}$ and
$\bfV_1 = \wh \bfV \bfQ_1 \in \R^{n\times r_1}$.
	\end{enumerate} 
	
	The time-updated numerical solution is then obtained as $\bfY_1 = \bfU_1\bfS_1\bfV_1^{\top}$.
	
\subsection{Rationale for choosing the augmented matrix \eqref{wh-S1}} 
\label{subsec:rationale}
In the rank-adaptive BUG integrator described in Section~\ref{sec:rabug}, a Galerkin method is done in the space spanned by the augmented left and right bases $\widehat \bfU$ and $\widehat \bfV$.  When these bases are chosen with the first $r$ columns formed by $\bfU_0$ and $\bfV_0$, respectively (as is done above), then the differential equation \eqref{S-step-rabug} for $\widehat\bfS(t)$ 
becomes
\begin{align}\label{eq:augmentedSMatrix}
    \dot{\widehat \bfS}(t) = \widehat \bfU^\top \bfF(t,\widehat \bfU \widehat \bfS(t) \widehat \bfV^\top) \widehat \bfV =
\begin{pmatrix} 
\bfU_0^\top \widehat \bfF(t) \bfV_0 & \bfU_0^\top \widehat \bfF(t) \widetilde \bfV_1
\\
\widetilde \bfU_1^\top \widehat \bfF(t) \bfV_0 & \widetilde \bfU_1^\top \widehat \bfF(t) \widetilde \bfV_1^\top
\end{pmatrix}
\end{align}
with $\widehat \bfF(t)=\bfF(t,\wh\bfU \wh\bfS(t) \wh \bfV^\top)$ for short, 
with the initial value
$$
\widehat \bfS(t_0) = 
\begin{pmatrix}
\bfS_0 & 0
\\
0 & 0
\end{pmatrix}.
$$
For the upper left block, this gives us, with $\bar \bfF(t)=\bfF(t,\bar \bfY(t))$ for $\bar \bfY(t)= \bfU_0 \bar \bfS(t) \bfV_0^\top$ and
up to an $O(h^2)$ approximation error, under the assumption of a bounded and Lipschitz-continuous $\bfF$ as will be made in the following,
$$
\widehat \bfS_{(0,0)}(t_1) = \bfS_0 + \int_{t_0}^{t_1} \bfU_0^\top \widehat \bfF(t) \bfV_0 \,dt \approx 
\bfS_0 + \int_{t_0}^{t_1} \bfU_0^\top \bar \bfF(t) \bfV_0 \, dt = \bar \bfS(t_1).
$$
For the lower left block
we have with $\bfF^K(t)= \bfF(t,\bfK(t)\bfV_0^\top)$ and again
up to an $O(h^2)$ approximation error, 
\begin{align*}
\widehat \bfS_{(1,0)}(t_1) &= \int_{t_0}^{t_1} \widetilde \bfU_1^\top \widehat \bfF(t) \bfV_0 \, dt
\\
&\approx  \int_{t_0}^{t_1} \widetilde \bfU_1^\top \bfF^K(t) \bfV_0 \, dt 
= \int_{t_0}^{t_1} \widetilde \bfU_1^\top \dot \bfK(t) \, dt = \widetilde \bfU_1^\top \bfK(t)= \bfS_1^K.
\end{align*}
Similarly we have up to $O(h^2)$
$$
\widehat \bfS_{(0,1)}(t_1) \approx \widetilde \bfS_1^L.
$$
Finally, the lower right block 
\begin{equation}\label{S-11}
\widehat \bfS_{(1,1)}(t_1) = \int_{t_0}^{t_1} \widetilde \bfU_1^\top \widehat \bfF(t) \widetilde \bfV_1 \, dt
\end{equation}
is of size $O(h\eps)$ if $\widehat\bfF(t_0)=\bfF(t_0,\bfY_0)$ is tangential up to $O(\eps)$ (with a small $\eps>0$) to the rank-$r$ matrix manifold at $\bfY_0=\bfU_0\bfS_0\bfV_0$. This holds true because for tangential matrices $\bfT$ at $\bfY_0$, which are of the form
$\bfT=\bfK\bfV_0^\top+\bfU_0\bfL^\top$ with arbitrary matrices $\bfK$ and $\bfL$ having $r$ columns, we have 
$\widetilde \bfU_1^\top \bfT \widetilde \bfV_1^\top = \bfzero$ as a consequence of the orthogonality relations $\widetilde \bfU_1^\top \bfU_0=\bfzero$ and $\widetilde \bfV_1^\top \bfV_0=\bfzero$. Approximate tangentiality to a low-rank matrix manifold of the vector field $\bfF(t,\bfY)$ near the solution of \eqref{ode} is a condition that needs to be imposed for dynamical low-rank approximation, since the given vector field is projected onto the tangent space at the current approximation, at every instance of time.

\begin{remark}
An alternative rank-adaptive algorithm, which is closely related to the rank-adaptive BUG integrator, can be constructed by not setting the term $\widehat \bfS_{11}(t_1)$ to zero, but instead computing it in a serial step according to
\begin{align}
   \dot{\widehat \bfS}_{(1,1)} = \widetilde\bfU_1^\top \widehat \bfF(t) \widetilde \bfV_1,\quad  \widehat \bfS_{(1,1)} = \mathbf{0}\in\mathbb{R}^{r\times r},
\end{align}
or simply setting
\begin{align}
\label{S-simple}
\widehat \bfS_{11}(t_1) = h\, \widetilde\bfU_1^\top \bfF(t_0,\bfY_0) \widetilde \bfV_1.
\end{align}
The core difference to the rank-adaptive BUG integrator is that this method does not require a $2r \times 2r$ coefficient update, but instead performs two $r \times r$ Galerkin steps (at least approximately). As we will see in the next subsection, the term in \eqref{S-simple} is anyway computed to decide on whether the rank needs to be increased.
\end{remark}

\subsection{Step rejection strategy}
\label{subsec:rejection}
To increase the robustness of the method, we propose to repeat the step in the following cases. This equally applies to the rank-adaptive BUG integrator of Section~\ref{sec:rabug}.
\begin{enumerate}
\item If $r_1= 2r$, then repeat the step with the augmented left and right bases $\widehat \bfU$ and $\widehat \bfV$ in the role of $\bfU_0$ and $\bfV_0$, respectively.
\item If $r_1 < 2r$, compute $\eta= \| \widetilde \bfU_1^\top \bfF(t_0,\bfY_0) \widetilde \bfV_1^\top\|$ and if $h\eta > c\vartheta$ (e.g., $c=10$) with the tolerance $\vartheta$ used in \eqref{vartheta}, then repeat the step with the augmented left and right bases $\widehat \bfU$ and $\widehat \bfV$ in the role of $\bfU_0$ and $\bfV_0$.
\end{enumerate}

{\it ad} 1.
The proposed rank-adaptivity works remarkably well in many cases, as reported in \cite{ceruti2022rank} for the integrator of Section~\ref{sec:rabug}, but it has been noted in \cite{Hau22} that it does not work satisfactorily in cases where a steep increase of the rank is necessary within a step. This is because the rank cannot be more than doubled in a step, and this may not be sufficient, for example when starting from rank~1. A  simple remedy is to repeat the step when no singular values are truncated, which indicates that an even higher rank is required.
With this criterion, an arbitrary rank increase is possible from one step to the next, if needed.

{\it ad} 2. Good approximation can be expected only in the case of a small normal component of $\bfF(t,\bfY)$ with respect to the tangent space of the rank-$r$ matrix manifold  at $(t,\bfY)$; cf.~\cite{kieri2016robust}.
More precisely, with a small $\eps>0$ we need
\begin{equation}\label{eps}
\| \mathrm{P}_{r}^\perp(\bfY_0) \bfF(t_0,\bfY_0) \| \le \eps,
\end{equation}
where 
$\mathrm{P}_{r}^\perp(\bfY)=\text{I}-\mathrm{P}_{r}(\bfY)$ and $\mathrm{P}_{r}(\bfY)$ is the orthogonal projection onto the tangent space $T_\bfY \mathcal{M}_{r}$ 
of the rank-${r}$ matrix manifold $\mathcal{M}_{r}$ at $\bfY\in \mathcal{M}_{r}$. 

The algorithm as it stands does not control $\| \mathrm{P}_{r}^\perp(\bfY) \bfF(t,\bfY)\|$. While this quantity is not computationally accessible with a computational cost that is linear in $m+n$, we can inexpensively compute the normal component of the orthogonal projection $\wh \bfU\wh \bfU^\top \bfF_0 \wh\bfV\wh\bfV^\top$ of $\bfF_0=\bfF(t_0,\bfY_0)$  onto the augmented subspace spanned by $\wh \bfU$ and $\wh\bfV$, which is shown to be (by the proof below)
\begin{equation} \label{P-perp-F}
\| \mathrm{P}_{r}^\perp(\bfY_0) [\wh \bfU\wh \bfU^\top \bfF_0 \wh\bfV\wh\bfV^\top] \| = \| \widetilde \bfU_1^\top \bfF_0 \widetilde \bfV_1\|=\eta.
\end{equation}
Taking $\wh \bfU\wh \bfU^\top \bfF_0 \wh\bfV\wh\bfV^\top$ as a reduced model of $\bfF_0$ in the computationally inaccessible term $\| \mathrm{P}_{r}^\perp(\bfY_0) \bfF_0\|$ and using \eqref{P-perp-F} leads us to the step rejection criterion 2. This criterion is even more important for the new parallel integrator than for the integrator of Section~\ref{sec:rabug} because the term $\wh S_{(1,1)}$ of \eqref{S-11}, which has norm $h\eta+O(h^2)$, is discarded in the algorithm.

Concerning the extra computational cost of computing $\eta$, we note that this is inexpensive if $\bfF_0$ is available (or is approximated) as a factorized low-rank matrix $\bfF_0=\bfG \bfH^\top$ with slim matrices $\bfG$ and $\bfH$, since then 
$\widetilde \bfU_1^\top \bfF_0 \widetilde \bfV_1=(\widetilde \bfU_1^\top \bfG) (\widetilde \bfV_1^\top \bfH)^\top$ is the product of two small matrices and the factors $\bfG$ and $\bfH$ need not be recomputed, but are still available from Step 1 of the integrator.  If instead only matrix-vector products $\bfF_0\bfv$ are accessible, then $\bfF_0 \widetilde \bfV_1$ must be computed anew. In this case, and unless the rank $r$ is very small, it can suffice to roughly approximate $\eta$ by taking only a few columns of $\widetilde \bfU_1$ and $\widetilde \bfV_1$ instead of all $r$ columns.

{\it Proof of \eqref{P-perp-F}.} It is known from \cite{koch2007dynamical} that for $\bfY_0=\bfU_0\bfS_0\bfV_0^\top\in\mathcal{M}_r$ 
and $\bfZ\in\R^{m\times n}$,
$$
\mathrm{P}_{r}^\perp(\bfY_0)\bfZ =(\bfI-\bfU_0\bfU_0^\top) \bfZ (\bfI-\bfV_0\bfV_0^\top) .
$$
Here, $\bfI-\bfU_0\bfU_0^\top$ is the orthogonal projection onto the orthogonal complement of the range of $\bfU_0$. 
We choose a basis $\bfU_0^\perp$ of this space whose first $r$ column vectors are those of $\widetilde\bfU_1$ and the remaining columns are orthogonal to $\wh\bfU=(\bfU_0,\widetilde\bfU_1)$:
$$
\bfU_0^\perp = (\widetilde\bfU_1, \wh\bfU^\perp) \quad\text{with}\quad (\wh\bfU^\perp)^\top \wh\bfU=0.
$$
We then have
$
\bfI-\bfU_0\bfU_0^\top = \bfU_0^\perp (\bfU_0^\perp)^\top
$
and hence 
$$
(\bfI-\bfU_0\bfU_0^\top){\wh \bfU}\wh \bfU^\top= (\widetilde\bfU_1, \wh\bfU^\perp)
\begin{pmatrix} \bfzero & \bfI \\ \bfzero & \bfzero \end{pmatrix} 
\begin{pmatrix} \bfU_0^\top \\ \widetilde\bfU_1^\top \end{pmatrix} = \widetilde\bfU_1\widetilde\bfU_1^\top.
$$
Analogously we have
$$
(\bfI-\bfV_0\bfV_0^\top)\wh\bfV\wh \bfV^\top = \widetilde\bfV_1\widetilde\bfV_1^\top.
$$
So we obtain
\begin{align*}
 \mathrm{P}_{r}^\perp(\bfY_0) [\wh \bfU\wh \bfU^\top \bfF_0 \wh\bfV\wh\bfV^\top ]
 &= (\bfI-\bfU_0\bfU_0^\top)\wh\bfU\wh \bfU^\top \bfF_0 \bigl((\bfI-\bfV_0\bfV_0^\top)\wh\bfV\wh \bfV^\top\bigr)^\top
 \\
 &= \widetilde\bfU_1\widetilde\bfU_1^\top \bfF_0 \widetilde\bfV_1\widetilde\bfV_1^\top,
 \end{align*}
 which implies \eqref{P-perp-F}. \qed


\section{Robust error bound}
For the rank-adaptive BUG integrator described in Section~\ref{sec:rabug}, we derived in 
 \cite{ceruti2022rank} local and global error bounds that are robust to small singular values, based on corresponding results in
\cite{kieri2016robust} and \cite{ceruti2022unconventional}. With the interpretation of the parallel integrator as a perturbed BUG integrator given in Section~\ref{subsec:rationale}, 
such robust error bounds are also obtained for the parallel integrator.

\subsection{Local error bound}  Let $\wh\bfY_1^\mathrm{BUG}$ be the rank-augmented numerical result of the BUG integrator of Section~\ref{sec:rabug} after one step starting from $\bfY_0$, and as above $\wh\bfY_1$ that of the new parallel integrator. Assuming that $\bfF$ is bounded with bounded derivatives, Section~\ref{subsec:rationale} shows that with the bound \eqref{eps}, the difference between these two numerical results is bounded by
$$
\| \wh\bfY_1 - \wh\bfY_1^\mathrm{BUG} \| \le \wh c_1 h^2 + \wh c_2 h\eps,
$$
where $\wh c_1$ and $\wh c_2$ are independent of any singular values. It is readily verified that these constants actually only depend on a bound and a Lipschitz constant of $\bfF$. The local error bound for the rank-augmented BUG integrator given in  \cite{ceruti2022rank} is of the same type if the initial value $\bfA(t_0)$ is of rank $r$ and $\bfY_0=\bfA(t_0)$:
$$
\| \wh\bfY_1^\mathrm{BUG} - \bfA(t_1) \| \le  c_1^\mathrm{BUG} h^2 +  c_2^\mathrm{BUG} h\eps,
$$
where again $c_1$ and $c_2$ only depend on a bound and a Lipschitz constant of $\bfF$.
Furthermore, the rank-truncated matrix $\bfY_1$ computed in the algorithm with truncation tolerance $\vartheta$ satisfies 
$$
\| \bfY_1 - \widehat\bfY_1\| \le \vartheta.
$$
With the triangle inequality we thus obtain a robust local error bound for the parallel integrator:
\begin{equation}\label{local-error-bound}
\| \bfY_1 - \bfA(t_1)\| \le c_1 h^2 + c_2 h\eps + \vartheta,
\end{equation}
where $c_1$ and $c_2$ only depend on a bound and a Lipschitz constant of $\bfF$.

\subsection{Global error bound} From the local error bound we conclude to a bound of the global error $\bfY_n - \bfA(t_n)$ for $0\le nh\le T$ with a fixed final time $T$ using the familiar argument of
Lady Windermere's fan \cite[Section~II.3]{hairer1993} with error propagation by the exact flow. We do not give the details of the standard proof, but we
formulate the resulting global error bound and begin with the assumptions:
\begin{enumerate}
    \item $\bfF$ is bounded and Lipschitz continuous. 
   \item The normal component of $\bfF(t_n,\bfY_n)$ is $\eps$-small for $0\le nh \le T$: With the rank $r_n$ used in the $n$th time step,
   $$
   \| \mathrm{P}_{r_n}^\perp(\bfY_n) \bfF(t_n,\bfY_n) \| \le \eps.
   $$
   (Note that the step rejection strategy of Section~\ref{subsec:rejection} aims at ensuring $h\eps \le c \vartheta$.)
    \item The error of the starting value $\bfY_0$ of rank $r_0$ is $\delta$-small:
       $ \Vert \bfY_0 - \bfA(0) \Vert \leq \delta$.
\end{enumerate}
We then have the following global error bound.
\begin{theorem} Under assumptions $1.$~to $3.$, the error of the parallel rank-adaptive integrator with stepsize $h$ and rank-truncation tolerance $\vartheta$ is bounded by
$$
\| \bfY_n - \bfA(t_n)\|  \le C_0 \delta + C_1 h + C_2 \eps + C_3  \vartheta/h , \qquad 0\le nh \le T,
$$
where $C_i$ $(i=0,1,2,3)$ only depend on the bound and Lipschitz constant of $\bfF$ and on the length $T$ of the time interval.
\end{theorem}

\input{nums}


\section*{Conclusion} In the present contribution, a novel fully parallel numerical integrator for dynamical low-rank approximation has been presented. Through its construction as an approximation to the rank-adaptive BUG integrator, the new integrator inherits the desirable robust error bound of both the BUG and the projector-splitting integrators. 
The full parallelism comes at the cost of no longer fulfilling the exactness property, which however does not impact the solution quality in the investigated numerical experiments. Beyond the enhanced parallelism, there is also a computational advantage in the serial setting, since the proposed integrator does not need a coefficient update with the augmented left and right basis matrices.
The inherent parallelism of the new low-rank matrix integrator is expected to extend to yielding significant advantages for Tucker and tree tensor network formats, which we aim to explore in future work.


\bibliographystyle{siamplain}
\bibliography{references}
\end{document}

%% file: ex_shared.tex
\usepackage{lipsum}
\usepackage{amsfonts}
\usepackage{graphicx}
\usepackage{epstopdf}
\usepackage{algorithmic}
\ifpdf
  \DeclareGraphicsExtensions{.eps,.pdf,.png,.jpg}
\else
  \DeclareGraphicsExtensions{.eps}
\fi

\newsiamremark{remark}{Remark}
\newsiamremark{hypothesis}{Hypothesis}
\crefname{hypothesis}{Hypothesis}{Hypotheses}
\newsiamthm{claim}{Claim}

\headers{A parallel rank-adaptive integrator for DLRA}{}

\title{A parallel rank-adaptive integrator for dynamical low-rank approximation\thanks{Submitted to the editors April 12, 2023.
\funding{The work of J.~Kusch was funded by the Deutsche For\-schungs\-gemein\-schaft (DFG, German Research Foundation) - 491976834, and the work of Ch.~Lubich was supported by DFG via FOR 5413 - 465199066 and SFB-TRR 352 - 470903074. The work of G.~Ceruti was supported by the Swiss National Science Foundation (SNSF),
grant number 200020-178806.}}}

\author{Gianluca Ceruti\footnotemark[2]\thanks{Institute of Mathematics, EPF Lausanne, 1015 Lausanne, Switzerland.}
\and Jonas Kusch\footnotemark[3]\thanks{Numerical Analysis and Scientific Computing, Universit\"at Innsbruck, Technikerstraße 13, 6020 Innsbruck.}
\and Christian Lubich\footnotemark[4]\thanks{Mathematisches Institut, Universit\"at T\"ubingen, Auf der Morgenstelle 10, 72076 T\"ubingen, Germany.}
}

\usepackage{amsopn}
\DeclareMathOperator{\diag}{diag}

\def\R{{\mathbb R}}

\def\eps{\varepsilon}

\def\wh{\widehat}

\newcommand\bfv{{\mathbf v}}

\newcommand\bfI{{\mathbf I}}
\newcommand\bfA{{\mathbf A}}

\newcommand\bfD{{\mathbf D}}

\newcommand\bfF{{\mathbf F}}
\newcommand\bfG{{\mathbf G}}
\newcommand\bfH{{\mathbf H}}
\newcommand\bfK{{\mathbf K}}
\newcommand\bfL{{\mathbf L}}
\newcommand\bfM{{\mathbf M}}
\newcommand\bfN{{\mathbf N}}
\newcommand\bfP{{\mathbf P}}
\newcommand\bfQ{{\mathbf Q}}

\newcommand\bfS{{\mathbf S}}
\newcommand\bfT{{\mathbf T}}
\newcommand\bfU{{\mathbf U}}
\newcommand\bfV{{\mathbf V}}
\newcommand\bfW{{\mathbf W}}

\newcommand\bfY{{\mathbf Y}}
\newcommand\bfZ{{\mathbf Z}}

\newcommand\bfSigma{{\mathbf \Sigma}}

\newcommand\bfzero{0}

\def\eps{\varepsilon}

\def\phi{\varphi}

\def\diag{\mbox{diag}}



%% file: nums.tex
\section{Numerical results}
In the following, we provide numerical experiments to illustrate the behaviour of the new parallel integrator and compare it against the BUG integrator. The source code to reproduce all numerical experiments conducted in this manuscript is openly available and can be found in \cite{code}. 
\subsection{Planesource testcase of radiative transfer}
In the following, we investigate Ganapol's plane-source benchmark~\cite{ganapol2008analytical}. This test case considers an isotropic particle beam in slab geometry which is evolved in time by the radiative transfer equation. This equation models the transport and interaction of radiation particles with a background medium. For a particle density (or angular flux) $f=f(t,x,\mu)$, where $t\in\R_+$ denotes time, $x\in [a,b]\subset\R$ denotes space and $\mu\in[-1,1]$ denotes the direction of flight (or angle), the radiation transport equation with isotropic scattering reads 
\begin{equation}
	\label{eq:rte}
	\begin{aligned}
		&\partial_t f + \mu \partial_x f + f = \frac{1}{2} \int_{-1}^{1} f \,d\mu,
		\qquad 
		(x, \mu) \in [a,b] \times [-1,1] , \\[1mm]
		& f(t = 0) = \frac{1}{3\sqrt{2\pi}\cdot10^{-2}}\exp\Big(-\frac{x^2}{18\cdot 10^{-4}} \Big).
	\end{aligned}
\end{equation}
To discretize the solution $f$ in space and angle, a spatial grid with $N_x$ cells and midpoints $x_1\leq x_2 \leq \cdots \leq x_{N_x}$ with equidistant spacing $\Delta x$ is chosen. In angle, a modal P$_N$ approximation with $N$ expansion coefficients (also known as moments) with respect to the normalized Legendre polynomials $p_i(\mu)$ is used. Then, in cell $I_j = [x_j-\Delta x/2,x_j+\Delta x/2]$ the entries of the solution matrix $\mathbf{Y} \in\R^{N_x\times N}$ are given as
\begin{align}\label{eq:matrixODErad1D}
    Y_{ji}(t) \simeq \frac{1}{\Delta x}\int_{I_j}\int_{-1}^1 f(t,x,\mu) p_{i-1}(\mu)\,\mathrm{d}\mu\mathrm{d}x.
\end{align}
An upwind discretization in space yields the matrix ODE
\begin{align}\label{eq:Fkin}
    \dot{\bfY} = -\bfD_x \bfY \bfA^{\top} + \bfD_{xx} \bfY |\bfA|^{\top} + \bfY\bfG =: \mathbf{F}(\bfY) ,
\end{align}
where $\bfA\in\R^{N\times N}$ is the flux matrix with entries $a_{\ell k} = \int_{-1}^1 \mu p_{\ell}(\mu)p_{k}(\mu)\,\mathrm{d}\mu$ and $\bfG = \diag(0,1,\cdots,1)\in\R^{N\times N}$ is the scattering matrix. $\bfA$ is diagonalizable with real eigenvalues, i.e., $\bfA = \bfT\bm{\Lambda}\bfT^{-1}$ where $\bfT$ is a transformation matrix and $\bm{\Lambda}$ is a real-valued diagonal matrix. The absolute value of the flux matrix is then defined as $|\bfA| := \bfT|\bm{\Lambda}|\bfT^{-1}$. The tridiagonal stencil matrices are given as
\begin{align*}
D_{x,j,j\pm 1} = \pm\frac{1}{2\Delta x},\quad D_{xx,j,j\pm1} = \frac{\Delta x}{2h}, \quad\text{and }\enskip D_{xx,jj} = \frac{\Delta x}{h}.
\end{align*}
A forward Euler time discretization of \eqref{eq:matrixODErad1D} is stable under the time step restriction $h = \mathrm{CFL} \cdot \Delta x$, where $\mathrm{CFL}\in (0,1]$. It is straightforward to show that the parallel integrator will inherit this CFL conditon: Note that the proposed integrator performs a Galerkin update with the basis  
\begin{align*}
    \bfW = [\bfU_0 \otimes \bfV_0, \widetilde\bfU_1 \otimes \bfV_0, \bfU_0 \otimes \widetilde\bfV_1] \in \R^{N_x \cdot N \times 3r^2}.
\end{align*}
We vectorize the flux matrix $\mathbf{F}(\bfY)$ at $\bfY_0$ as $\mathbf{f}\in\R^{N_x\times N}$ according to
\begin{align*}
    f_{\mathrm{idx}(i,j)} = \left(\bfF(\bfY_0)\right)_{ij}, \quad \text{where } \mathrm{idx}(i,j) = (i-1)\cdot N_x + j
\end{align*}
as well as the numerical solution $\bfY_0$ as $\mathbf{y}\in\R^{N_x\times N}$ with entries $y_{\mathrm{idx}(i,j)} = \left(\bfY_0\right)_{ij}$. Then, with the Euclidean norm $\Vert \cdot \Vert_E$, we have
\begin{align*}
    \Vert \widehat\bfY_1 \Vert = \Vert \left(\mathbf{y} + h \mathbf{f}\right)\bfW \Vert_E \leq \Vert \mathbf{y} + h \mathbf{f} \Vert_E = \Vert \bfY_0 + h \bfF(\bfY_0) \Vert.
\end{align*}
From the time step restriction of the full problem we have $\Vert \bfY_0 + h \bfF(\bfY_0) \Vert \leq \Vert \bfY_0 \Vert$ and the truncation fulfills $\Vert \bfY_1 \Vert \leq \Vert \widehat\bfY_1 \Vert$. Hence $\Vert \bfY_1 \Vert \leq \Vert \bfY_0 \Vert$.


We are interested in the so-called scalar flux $\Phi(t,x) = \int_{-1}^{1} f(t,x,\mu) \,d\mu$ which is given on the discrete level as $\bm\Phi\in\R^{N_x}$ with entries $\Phi_j = Y_{j1}$. We choose a spatial domain of $[-5,5]$ which is discretized by $N_x=1000$ spatial cells. The angular dependence is represented by $N=500$ moments and the time step size is determined with $\mathrm{CFL} = 0.99$. As tolerance parameter we choose $\vartheta = \bar\vartheta\cdot \| \hat\bfSigma \|$, where $\bar\vartheta = 10^{-2}$. Moreover, the rejection step uses a value of $c = 1$ in the criterion $h\eta > c\vartheta$ and we use all $r$ basis vectors in $\widetilde U_1$ and $\widetilde V_1$ to compute $\eta$.

\begin{figure}[h!]
	\includegraphics[width=0.49\textwidth]{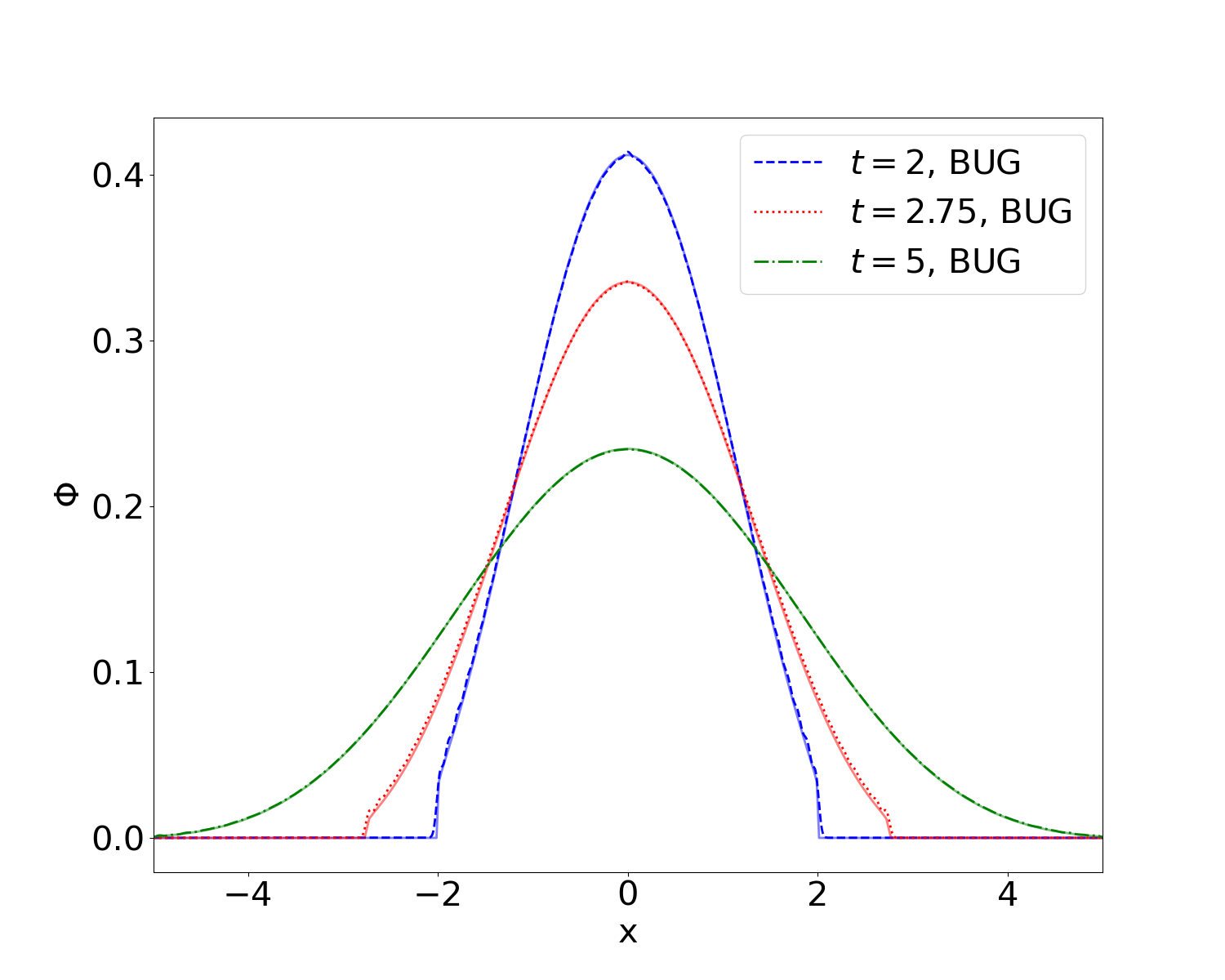}
	\includegraphics[width=0.49\textwidth]{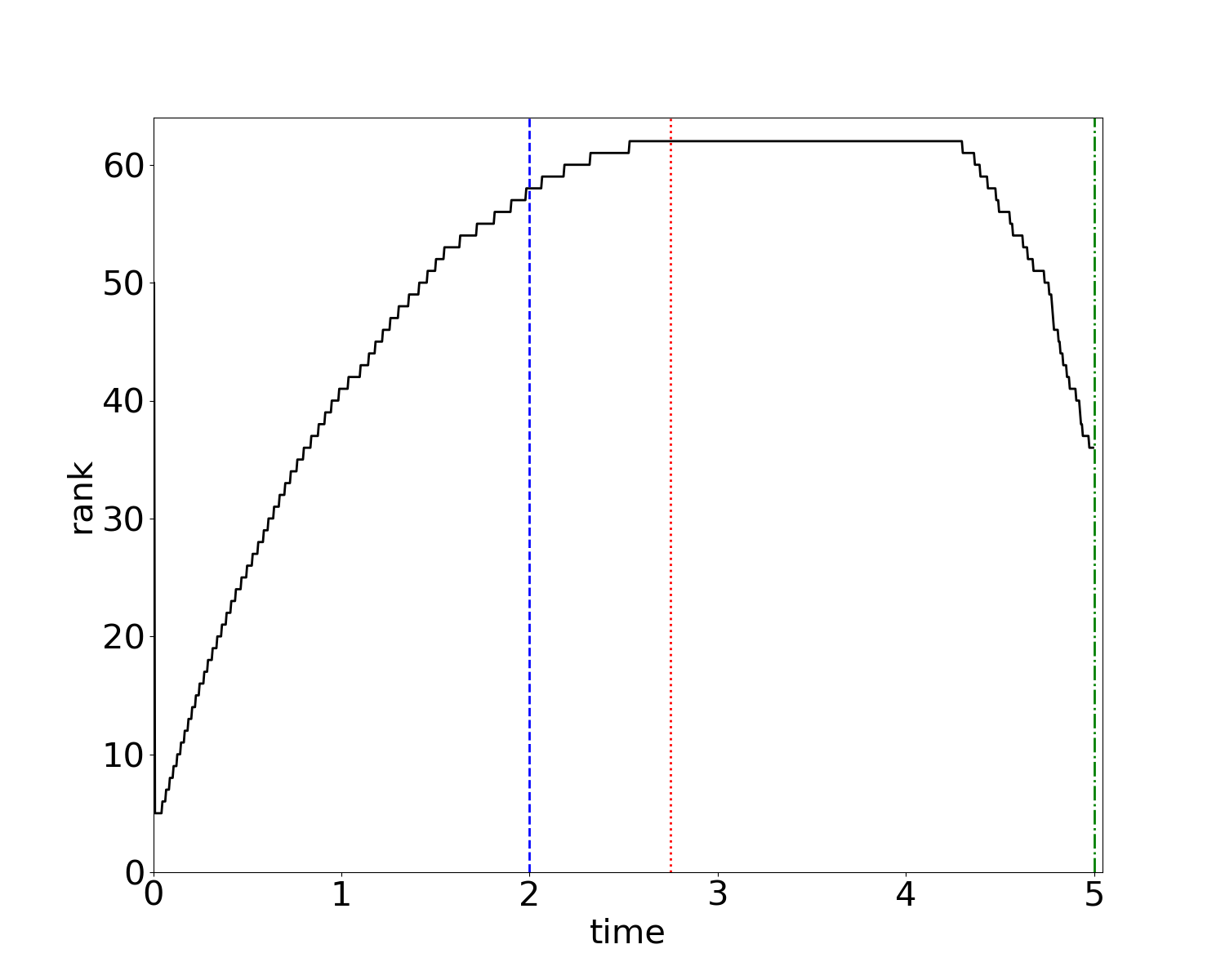}
		\includegraphics[width=0.49\textwidth]{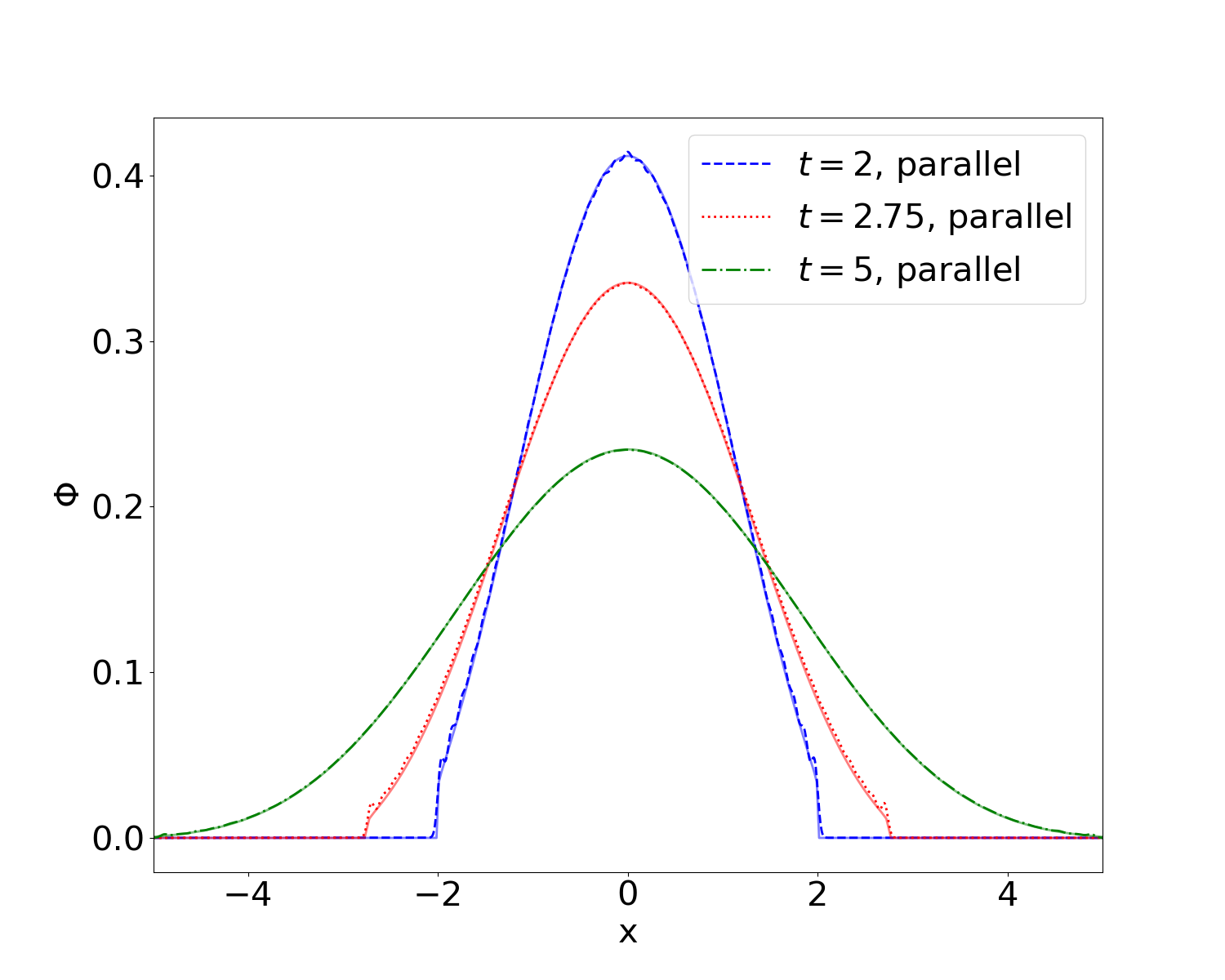}
	\includegraphics[width=0.49\textwidth]{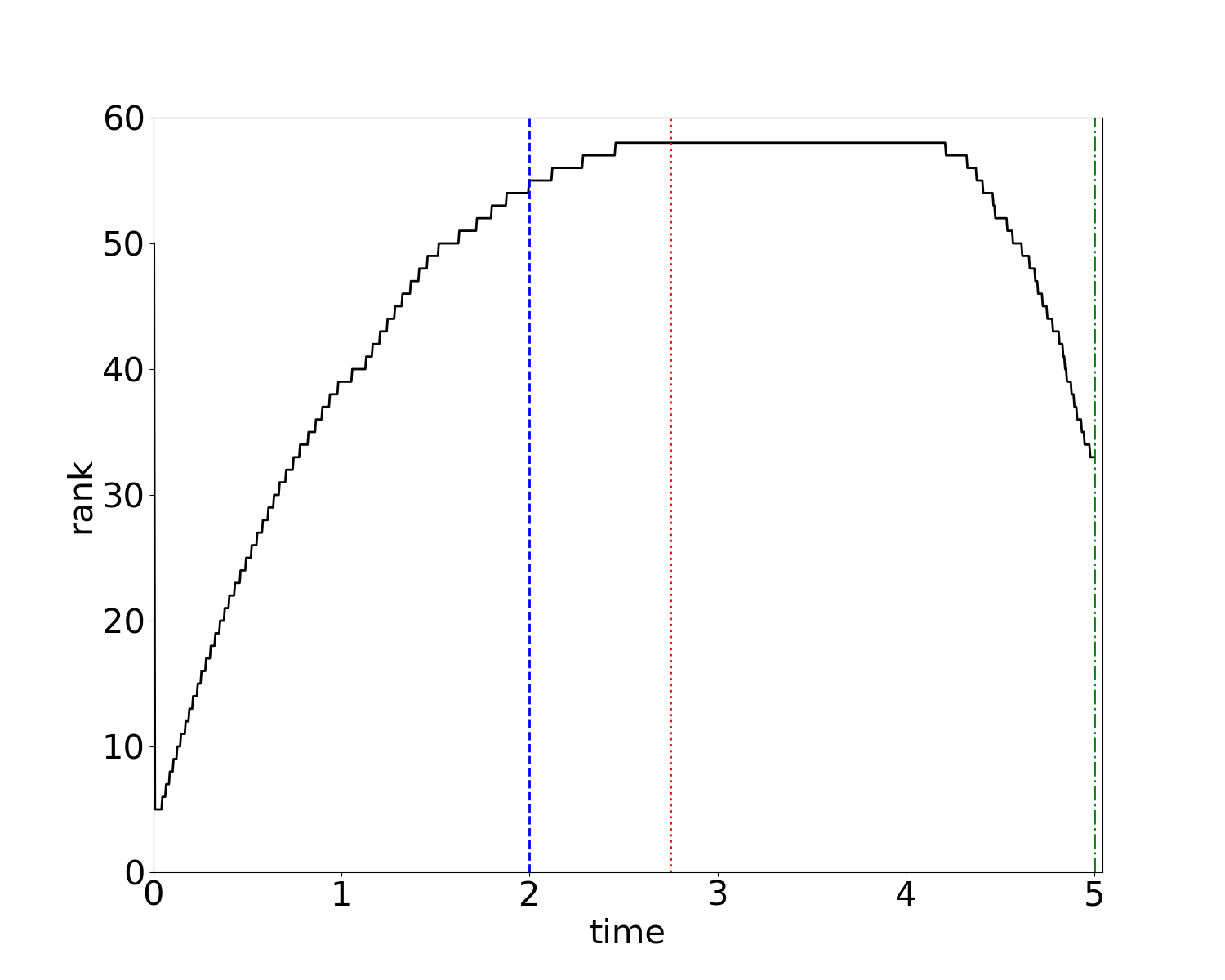}
	\caption{
		Top: Rank-adaptive BUG integrator. Bottom: Parallel integrator. Left: Scalar flux of the reference solution (solid lines) of the radiation transport equation \eqref{eq:rte} in comparison with the approximation of the DLRA integrators at different times. Right: Rank evolution of the approximation arising from the rank-adaptive integrator.}
	\label{fig:rte1Dphi}
\end{figure}

\begin{figure}[h!]
\centering
	\includegraphics[width=0.49\textwidth]{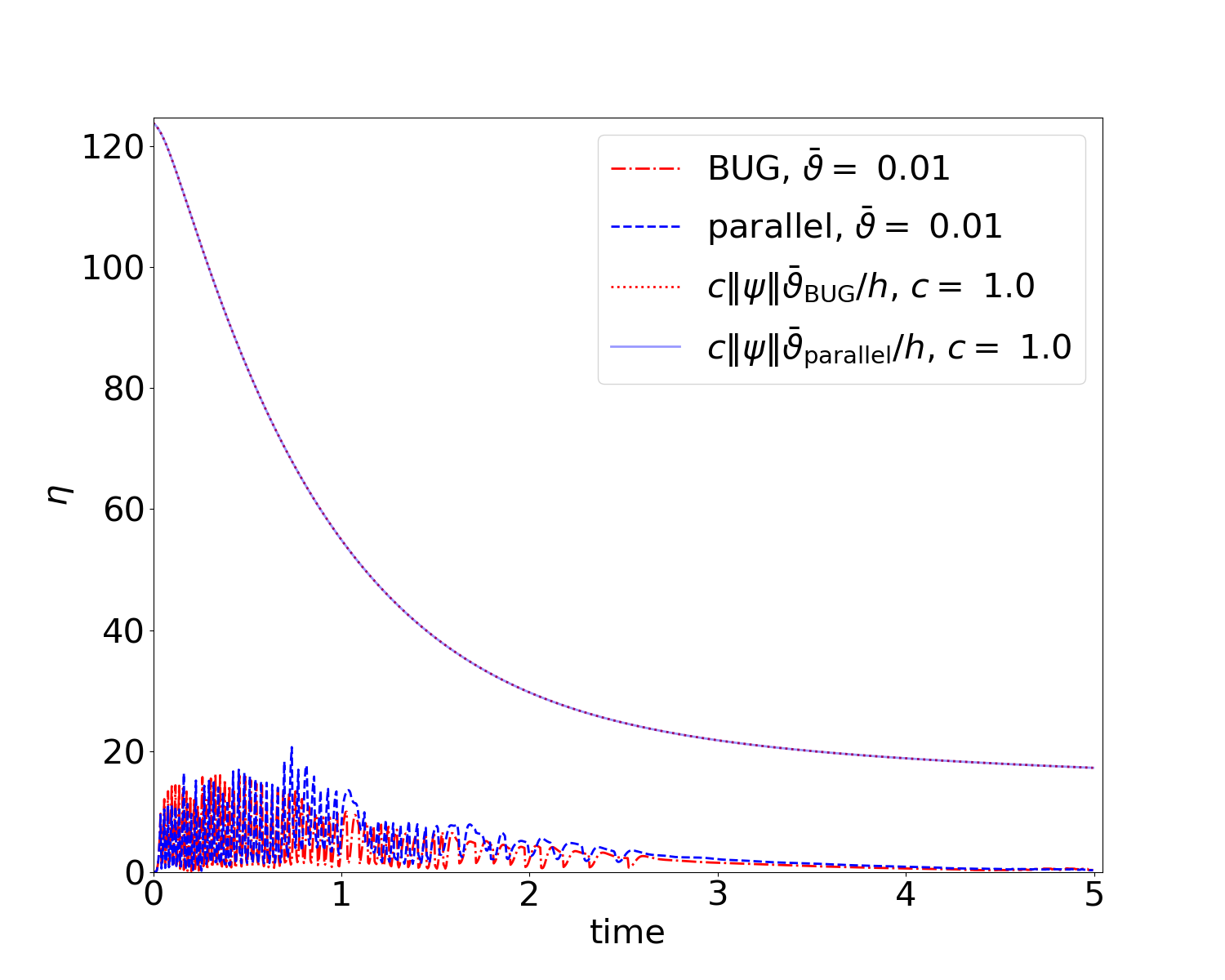}
	\caption{Recorded values for $\eta$ and rejection bound over time. It can be seen that even for a small value of $c=1$ the time update is never rejected.}
	\label{fig:rteEta}
\end{figure}

The left column of Figure~\ref{fig:rte1Dphi} depicts the scalar flux $\Phi$ at times $t \in \{2, 2.75, 5\}$ for the BUG and the new parallel integrator. The right column shows the rank evolution over time (top for BUG and bottom for the parallel integrator). It is observed that both integrators yield similar solution qualities as well as ranks. The resulting values for $\eta$ as well as the criterion of the rejection step $c\vartheta / h$ are presented in Figure~\ref{fig:rteEta}. We observe that the rejection step is not required in this testcase since the value of $\eta$ is sufficiently small throughout the entire time interval.
The parallel integrator at the coarse tolerance has a runtime of 26.8 seconds, whereas the rank-adaptive BUG integrator requires 47.6 seconds to run. At the fine tolerance, the parallel integrator computes the results in 39.4 seconds compared to 58.6 seconds for the rank-adaptive BUG integrator. That is, even at the fine tolerance, the parallel integrator is significantly faster than the rank-adaptive BUG integrator at the coarse tolerance.

\subsection{Linesource testcase of radiative transfer}
\begin{figure}
     \centering
     \begin{subfigure}[b]{0.49\textwidth}
         \centering
         \includegraphics[width=\textwidth]{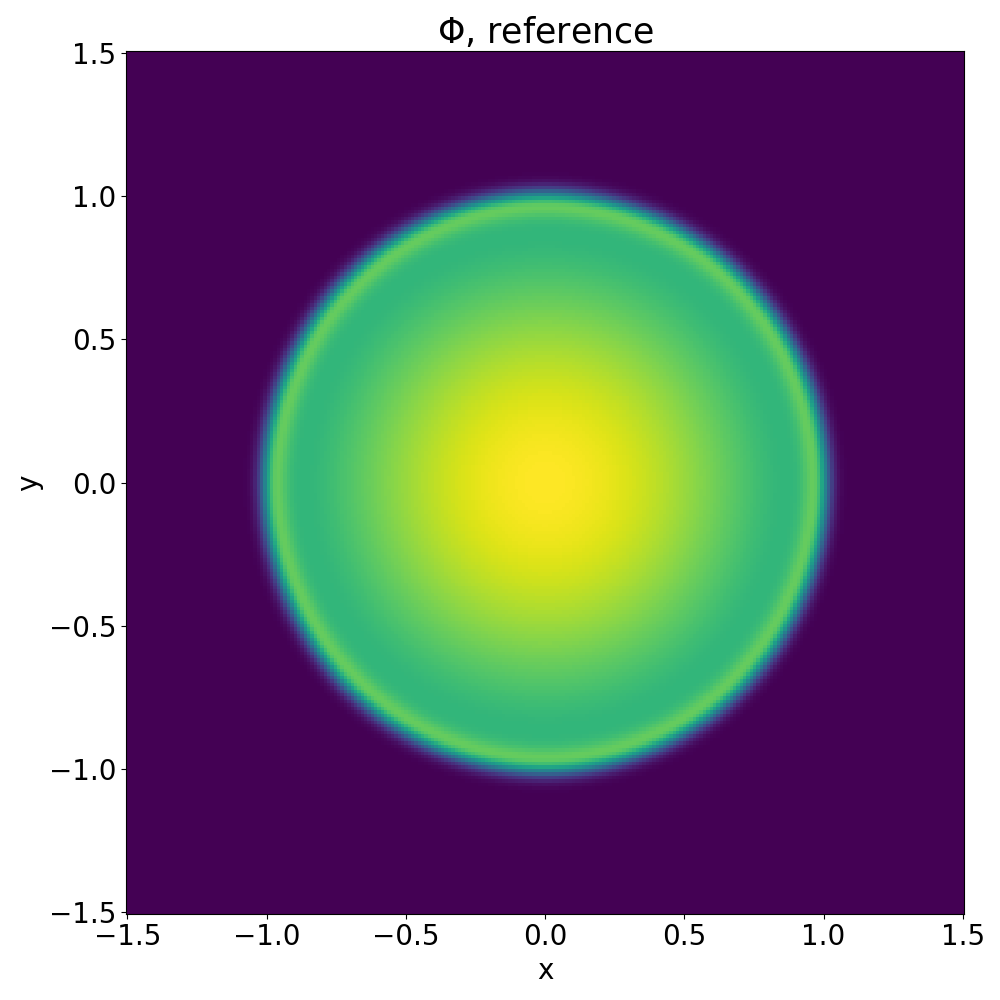}
         \caption{}\label{fig:scalar_flux_reference}
     \end{subfigure}
     \hfill
     \begin{subfigure}[b]{0.49\textwidth}
         \centering
         \includegraphics[width=\textwidth]{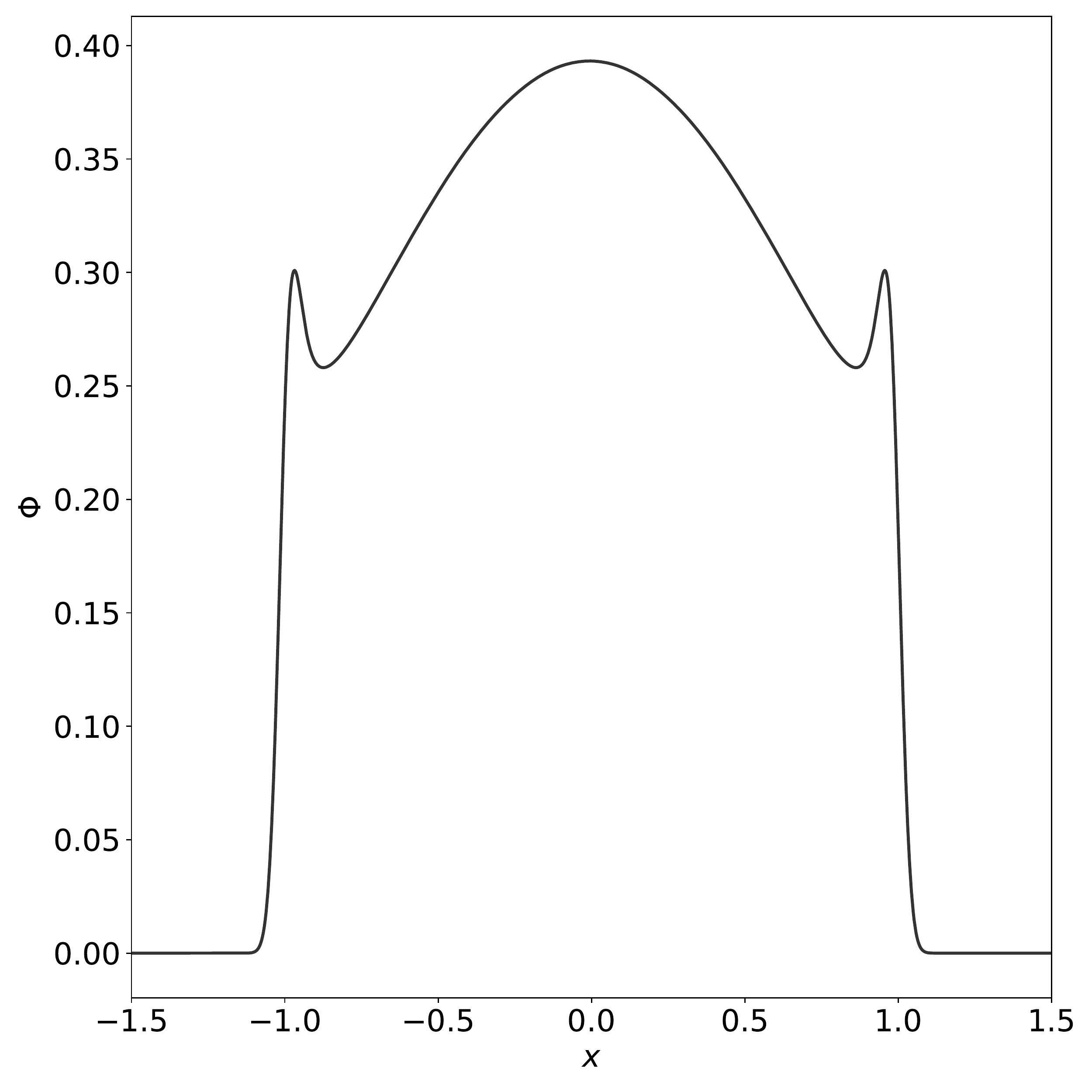}
         \caption{}
     \end{subfigure}
        \caption{
		Left: Scalar flux of reference solution. Right: Scalar flux of reference solution on a cut along $y=0$.}
	\label{fig:scalar_flux_linesource_ref}
\end{figure}

In the following, we extend the planesource testcase to a two-dimensional setting. This extension is known as the linesource testcase. For this benchmark, one approximates the solution $f = f(t,\mathbf{x},\mathbf{\Omega})$ to the two-dimensional radiative transfer equation
\begin{equation}
	\label{eq:rtels}
	\begin{aligned}
		&\partial_t f + \mathbf{\Omega}\cdot\nabla f + \sigma_t f = \frac{\sigma_s}{4\pi} \int_{\mathbb{S}^2} f \,d\mathbf{\Omega},\qquad 
		(\mathbf{x}, \mathbf{\Omega}) \in [-1.5,1.5]^2 \times \mathbb{S}^2 \;, \\[1mm]
		& f(t = 0) = \max\left\{10^{-4},\frac{1}{4\pi \sigma^2}\cdot \mathrm{exp}\left(-\frac{\Vert \mathbf{x}\Vert^2}{4\sigma^2}\right)\right\}.
	\end{aligned}
\end{equation}
Here, the two-dimensional spatial variable is denoted as $\mathbf x = (x,y)^{\top}$ and the angle (i.e., the direction of flight) is $\mathbf{\Omega}\in\mathbb{R}^3$. To model a constant solution along the $z$-axis, we set $\nabla = (\partial_x, \partial_y,0)^{\top}$. The initial condition is a Gaussian pulse with standard deviation $\sigma = 0.03$ and an isotropic distribution in direction of flight. The background medium through which particles move has a scattering cross-section $\sigma_s = 1$. The total cross section is defined as $\sigma_t = \sigma_s + \sigma_a$, where $\sigma_a = 0$ denotes the absorption cross-section. This benchmark is remarkably challenging for conventional numerical methods as it exposes deficiencies of individual discretizations. For more information, see e.g., \cite{garrett2013comparison}. This testcase is also challenging for low-rank methods as it requires high ranks to yield good approximations, see e.g., \cite{peng2020low,kusch2021robust,peng2023sweep}. At a final time $t_{\mathrm{end}}=1$, we are interested in the scalar flux $\Phi(t_{\mathrm{end}},\mathbf x) = \int_{\mathbb{S}^2} f(t,\mathbf x,\mathbf{\Omega}) \,\mathrm{d}\mathbf{\Omega}$. The linesource testcase is equipped with an analytic solution \cite{ganapol1999homogeneous,ganapol2008analytical} which is depicted in Figure~\ref{fig:scalar_flux_linesource_ref}.
To approximate the solution to the linesource testcase, we choose an equidistant spatial discretization of $250 \times 250$ spatial cells. The directional domain is discretized with a modal expansion in spherical harmonics basis functions. A polynomial degree of $39$ which corresponds to $1600$ expansion coefficients in angle is chosen. We use an upwind scheme with a forward Euler time integration to discretize space and time. As time step restriction, we use the CFL condition of the original problem and set $h = 0.7\cdot \Delta x$, where $\Delta x = \Delta y$ is the cell width in $x$. The chosen constant in the bound for $\eta$ is $c = 5$.

\begin{figure}
	\includegraphics[width=0.49\textwidth]{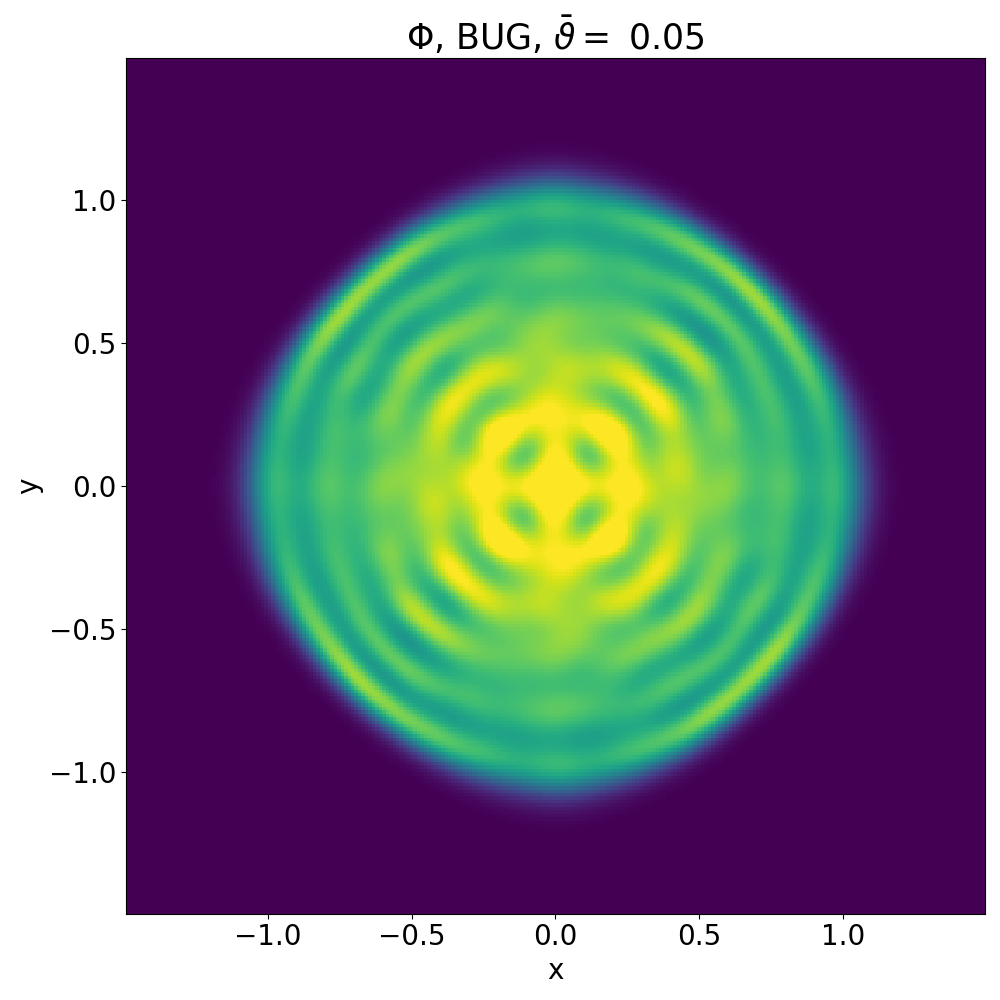}
    \includegraphics[width=0.49\textwidth]{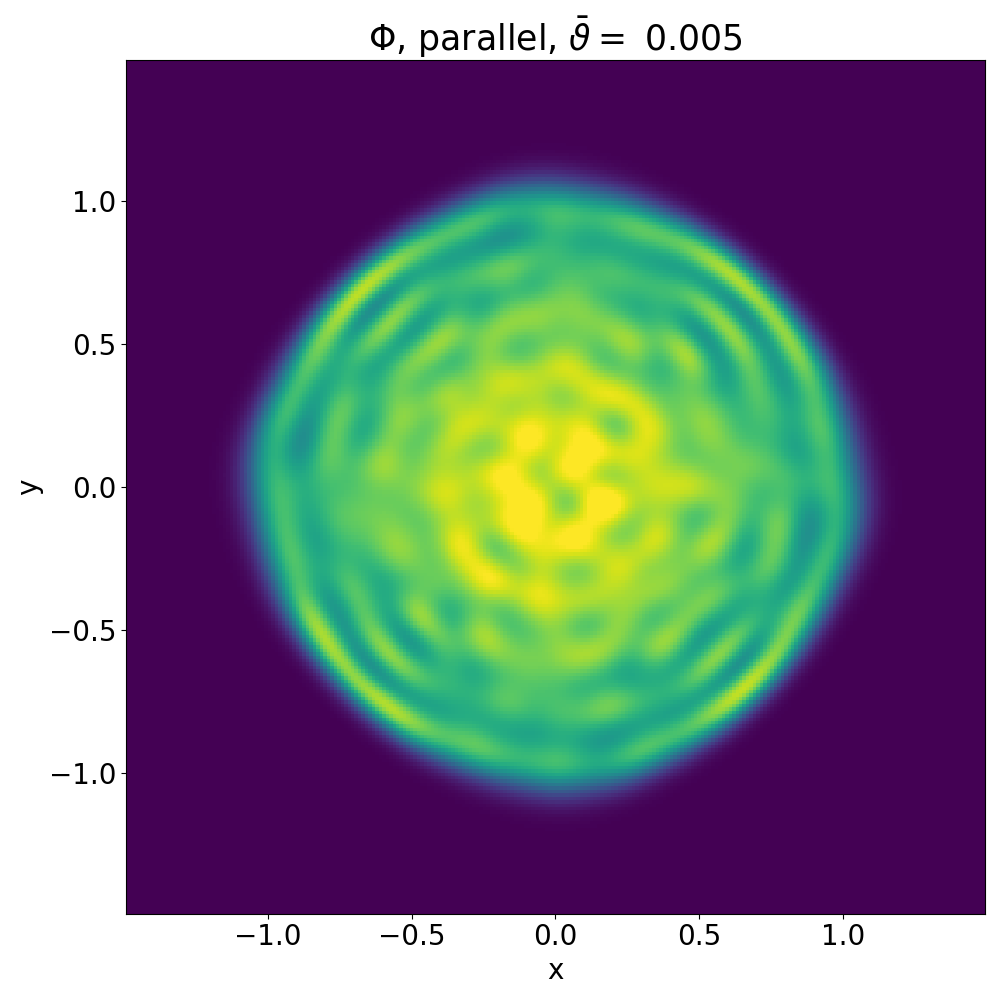}
	\includegraphics[width=0.49\textwidth]{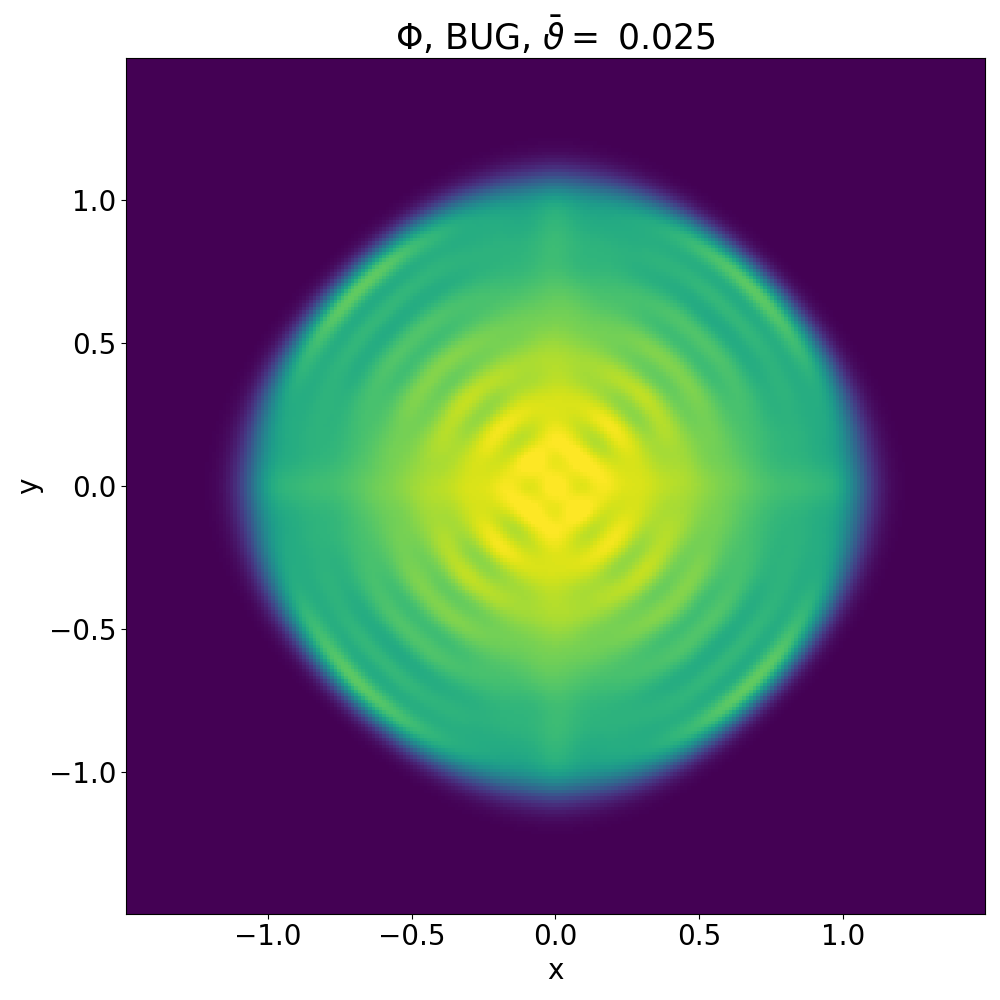}
    \includegraphics[width=0.49\textwidth]{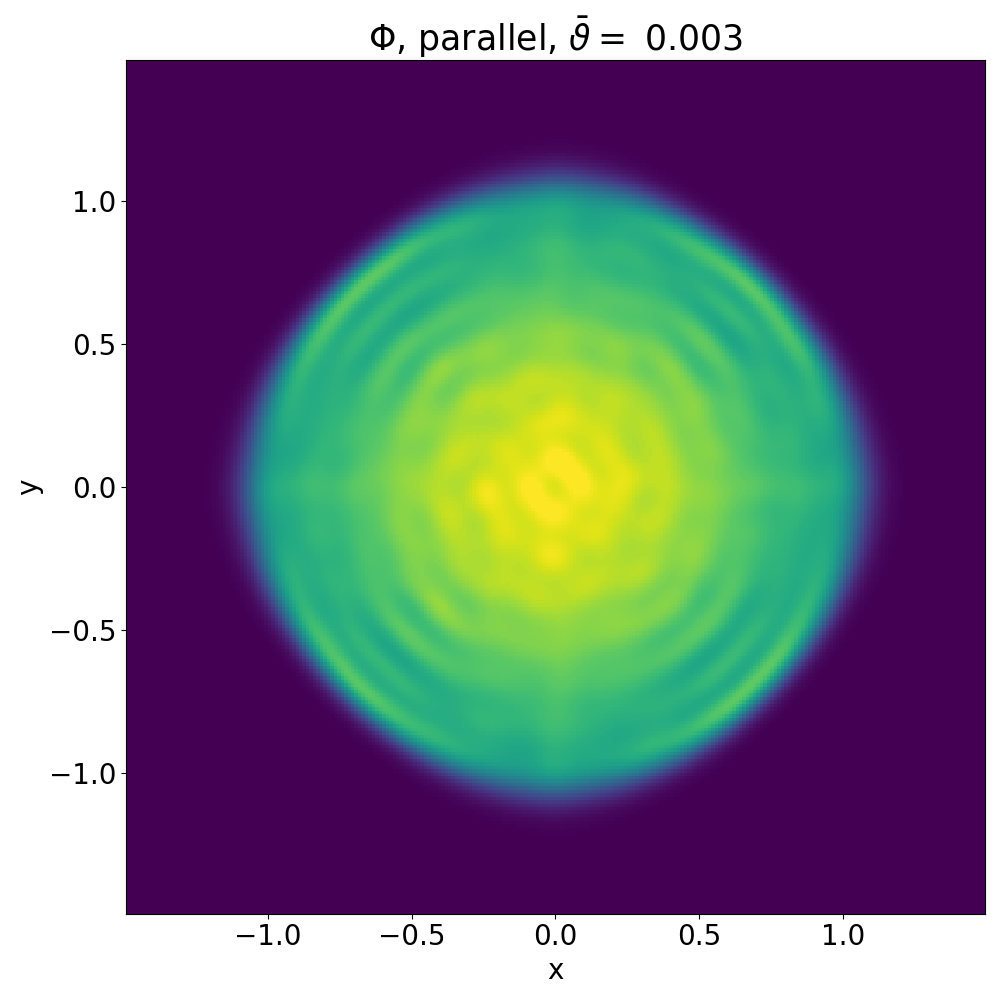}
    \includegraphics[width=0.49\textwidth]{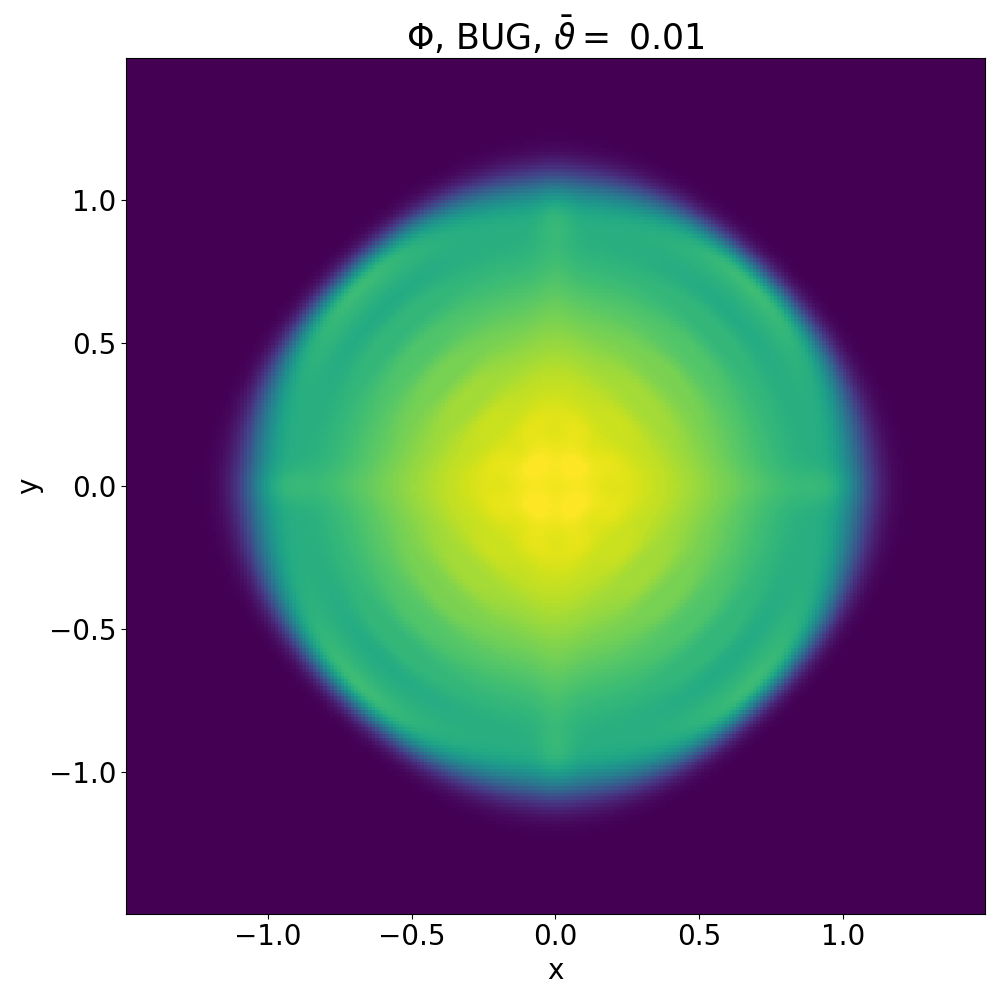}
    \includegraphics[width=0.49\textwidth]{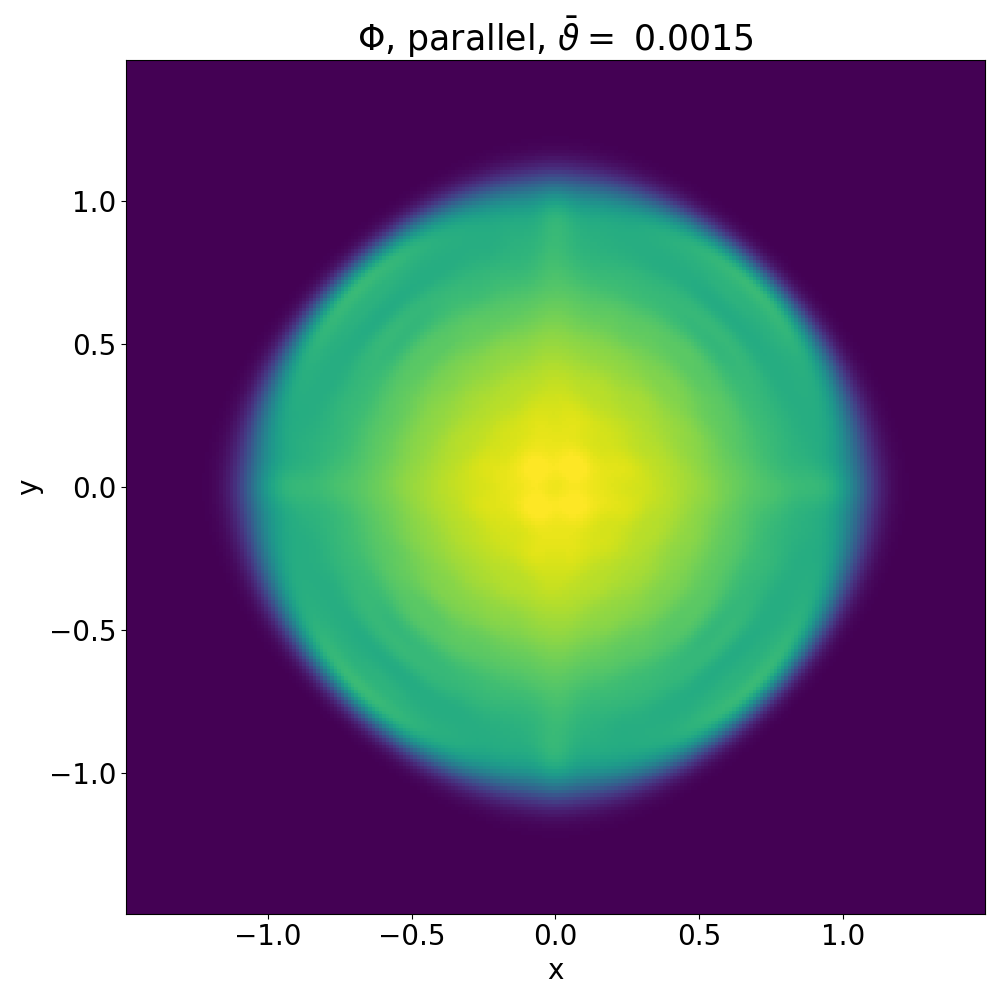}
	\caption{
		Left: Adaptive BUG integrator. Right: Parallel integrator.}
	\label{fig:linesource_DLRA}
\end{figure}

\begin{figure}
	\includegraphics[width=0.49\textwidth]{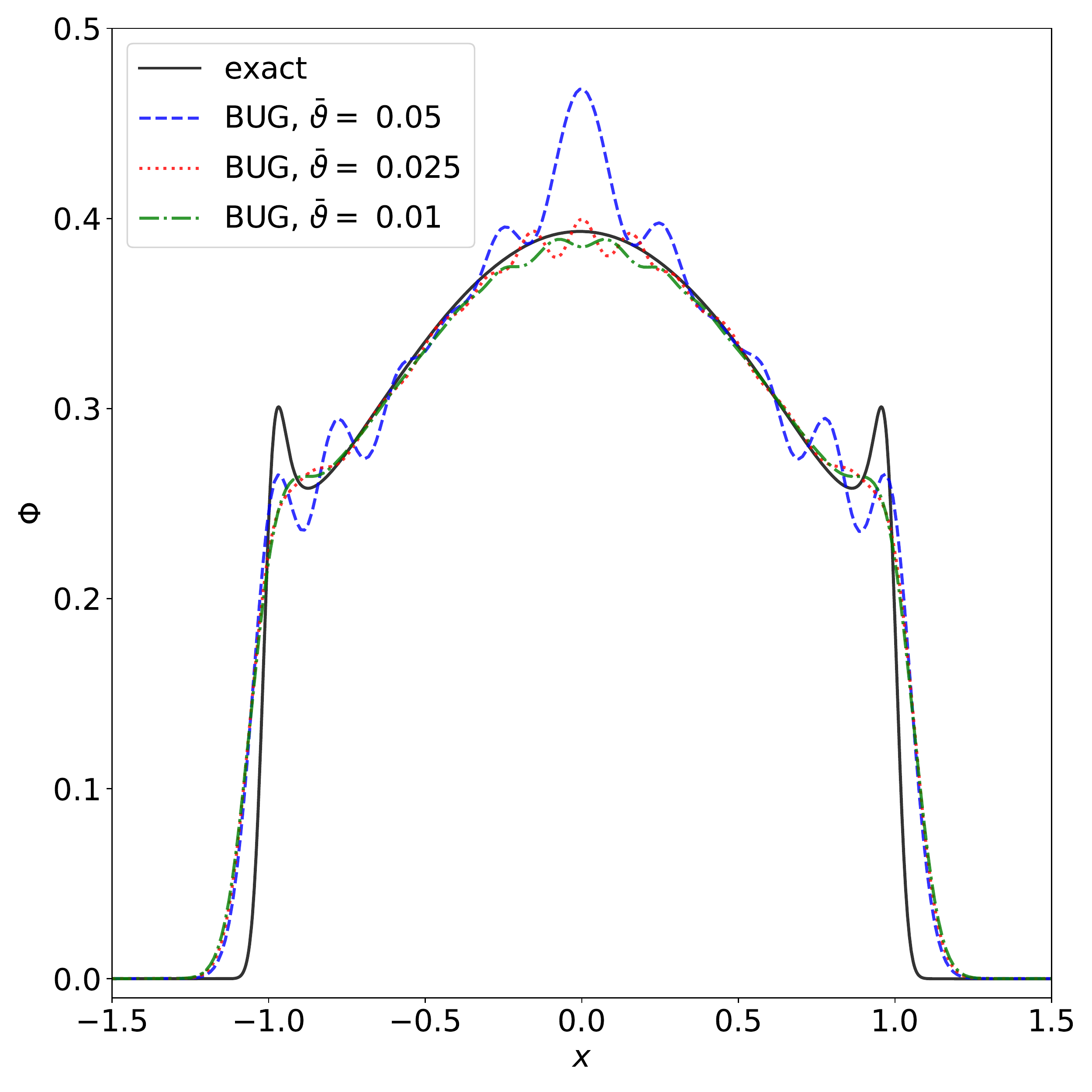}
    \includegraphics[width=0.49\textwidth]{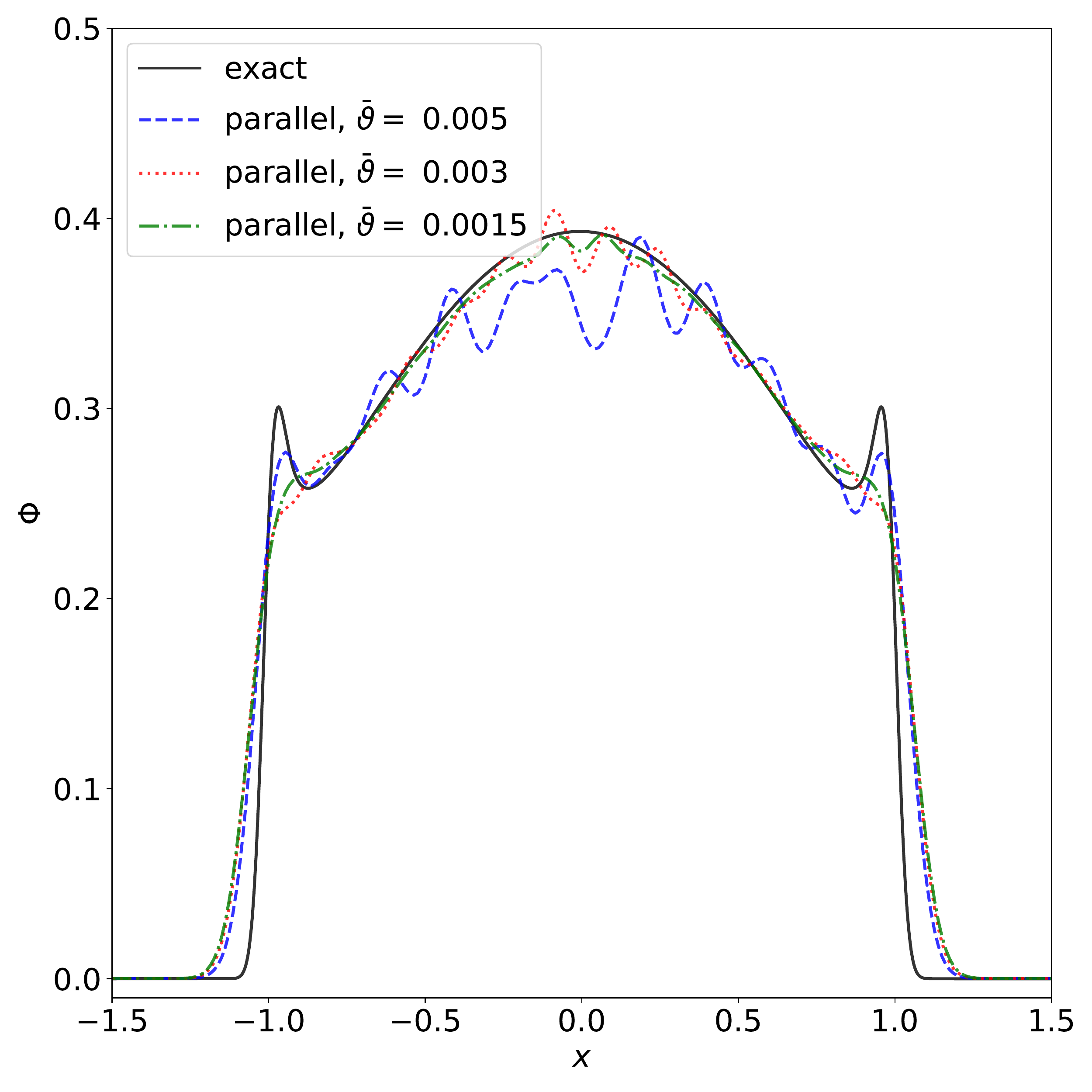}
	\caption{
		Left: Adaptive BUG integrator. Right: Parallel integrator.}
	\label{fig:linesource_DLRA_cut}
\end{figure}

\begin{figure}
\centering
	\includegraphics[width=0.49\textwidth]{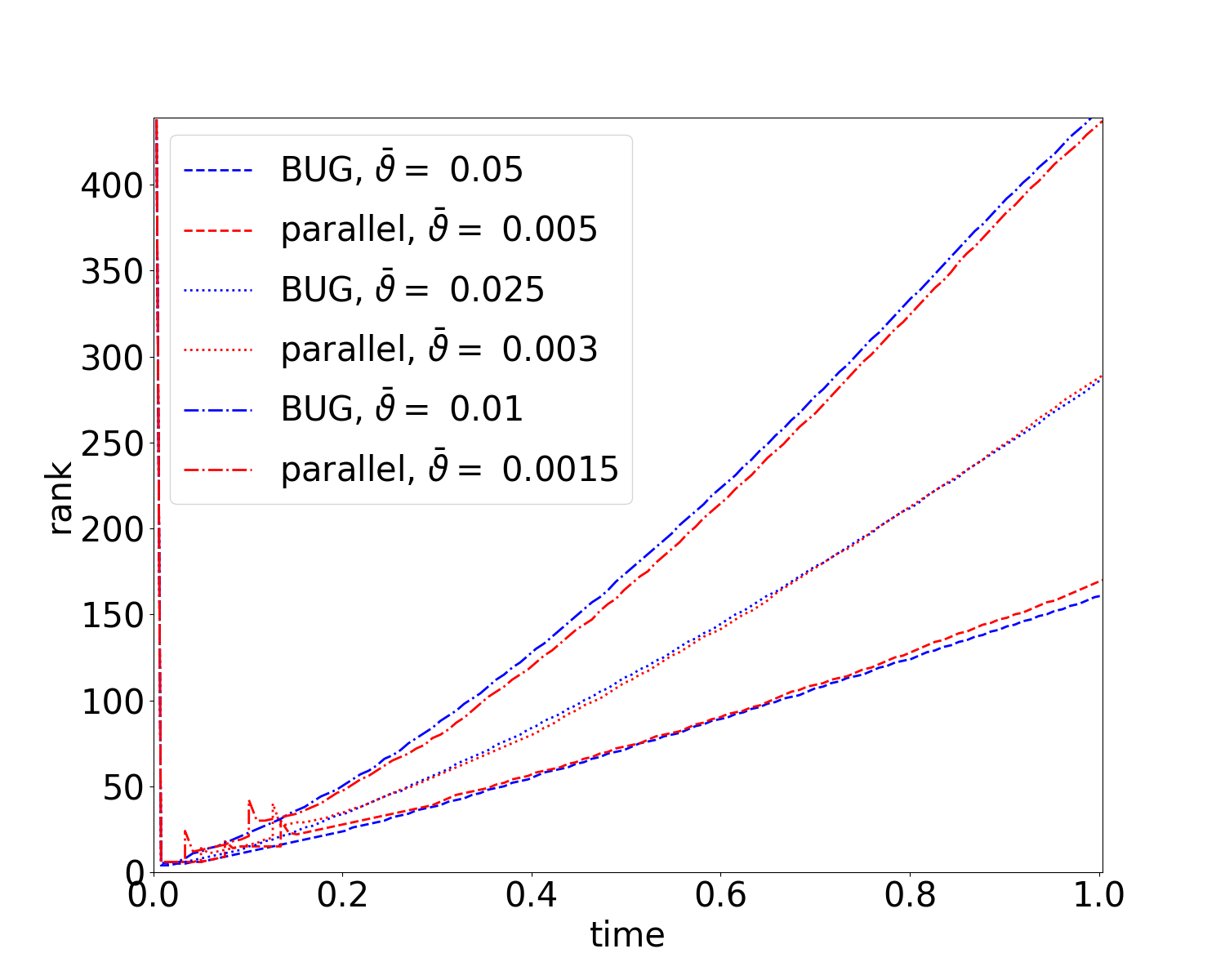}
 	\includegraphics[width=0.49\textwidth]{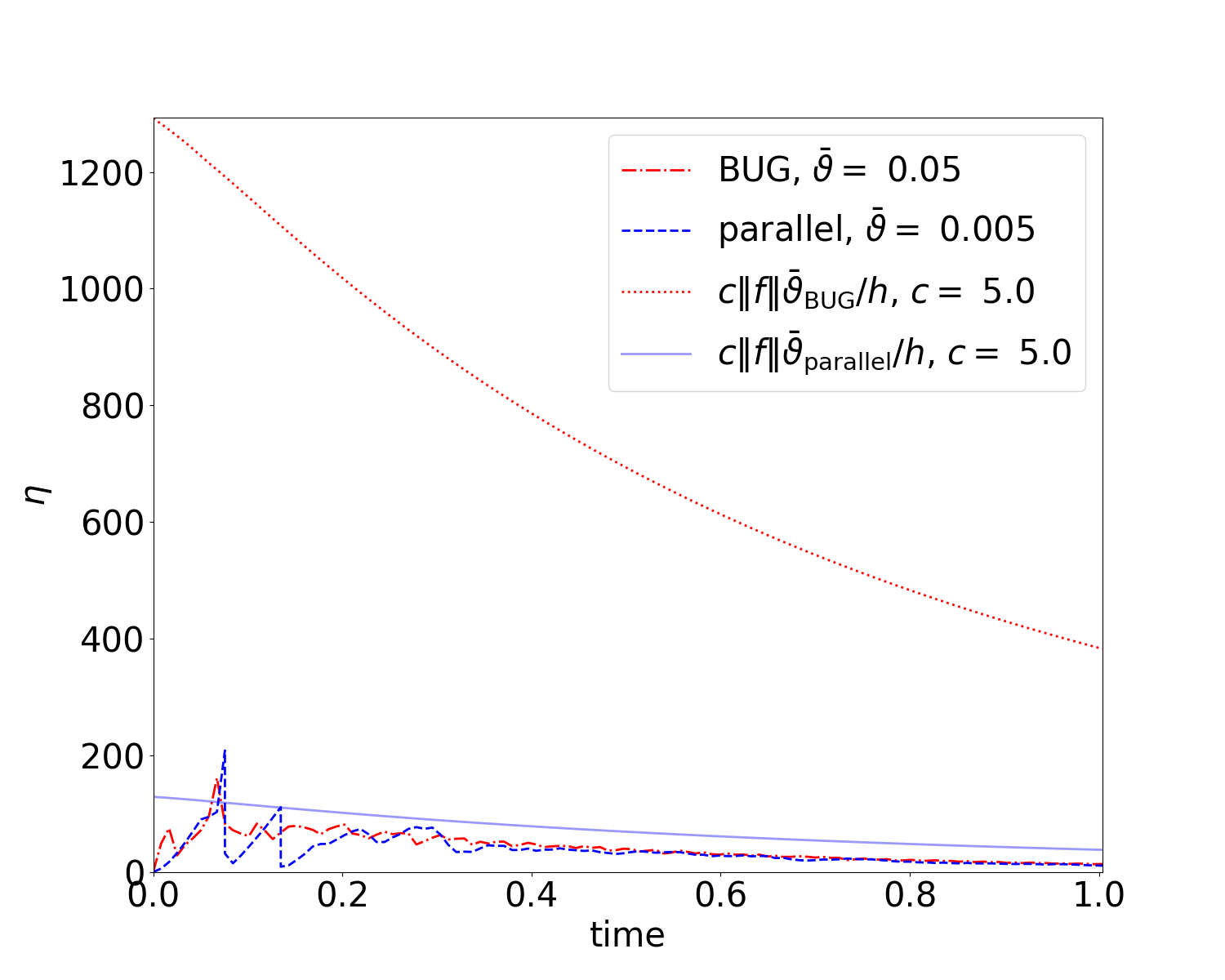}
  	\includegraphics[width=0.49\textwidth]{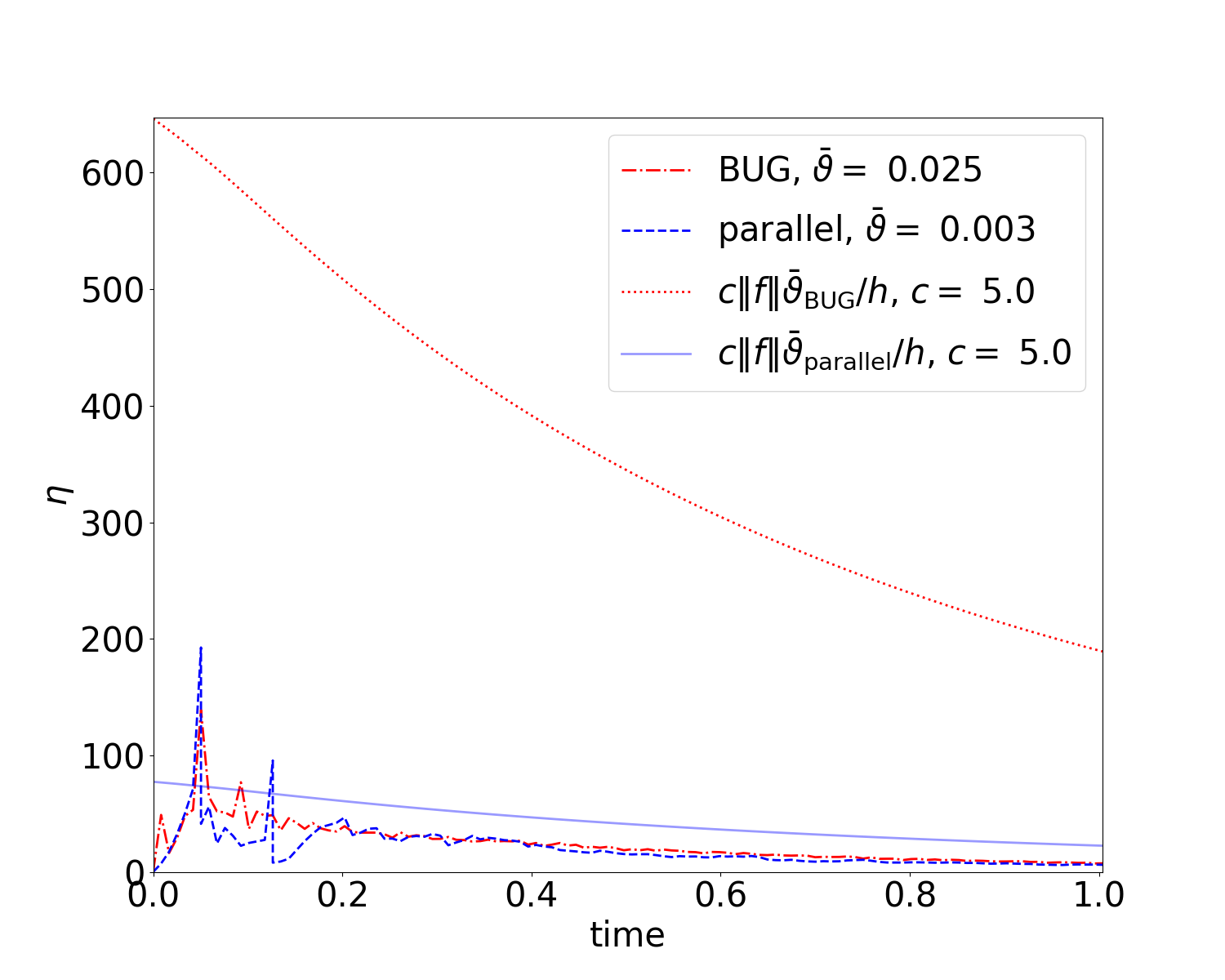}
   	\includegraphics[width=0.49\textwidth]{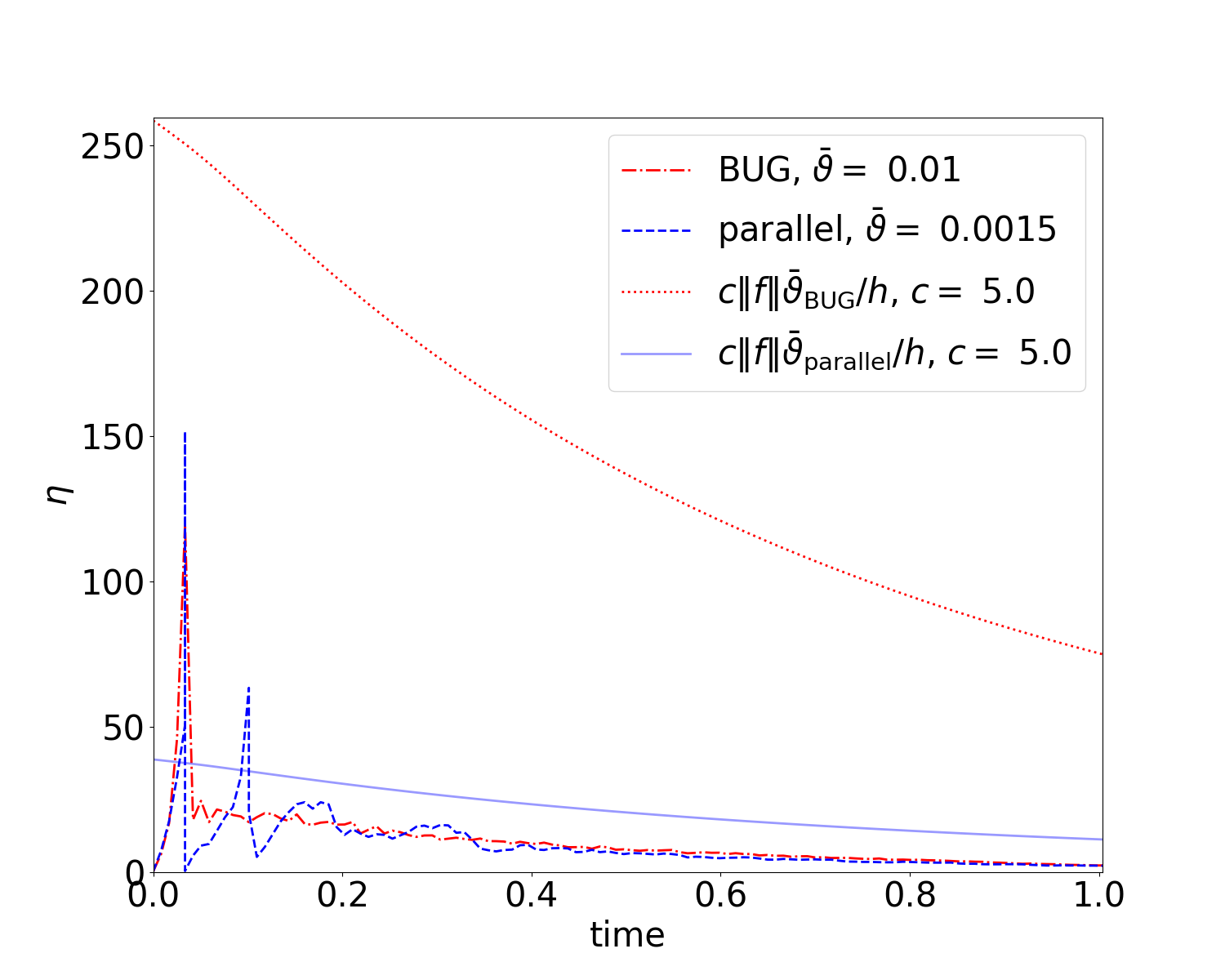}
	\caption{Rank and $\eta$ evolution in time.}
	\label{fig:rank_linesource}
\end{figure}

We observe that the ranks of the BUG as well as the parallel integrator will differ significantly when choosing the same tolerance parameter $\vartheta$ for both methods. Following \cite{ceruti2022rank}, we pick tolerances $\bar\vartheta\in\{0.05,0.025,0.01\}$ for the BUG integrator as well as $\bar\vartheta\in\{0.005,0.003,0.0015\}$ for the parallel integrator, which yields comparable ranks. The  resulting scalar fluxes are found in Figure~\ref{fig:linesource_DLRA}. Figure~\ref{fig:linesource_DLRA_cut} depicts the solution on a cut along $y=0$. The chosen ranks are depicted in Figure~\ref{fig:rank_linesource} together with the resulting time evolution of $\eta$.
We observe that both integrators yield similar numerical results at comparable ranks. Due to the fact that the parallel integrator does not necessitate a computation of a $2r \times 2r$ coefficient matrix update, its computational time is substantially decreased. At the coarse tolerance, the parallel integrator requires $186$ seconds compared to $256$ seconds for BUG. At the intermediate tolerance, the parallel integrator requires $308$ seconds compared to $462$ seconds for BUG. At the fine tolerance, the parallel integrator requires $506$ seconds compared to $782$ seconds for BUG. When the problem is sufficiently big and justifies the communication overhead of a parallel implementation, the use of the full parallelism will further benefit the parallel integrator. Such an implementation is however left for future work on parallel tensor integrators.

Note that the recorded values for $\eta$ depicted in Figure~\ref{fig:rank_lattice} indicate that a limited number of rejection steps are required for the parallel integrator in the beginning of the simulation. This is especially the case, since the chosen tolerance is significantly smaller compared to the tolerance of the BUG integrator.

\subsection{Lattice testcase of radiative transfer}
\begin{figure}
     \centering
     \begin{subfigure}[b]{0.49\textwidth}
         \centering
         \includegraphics[width=\textwidth]{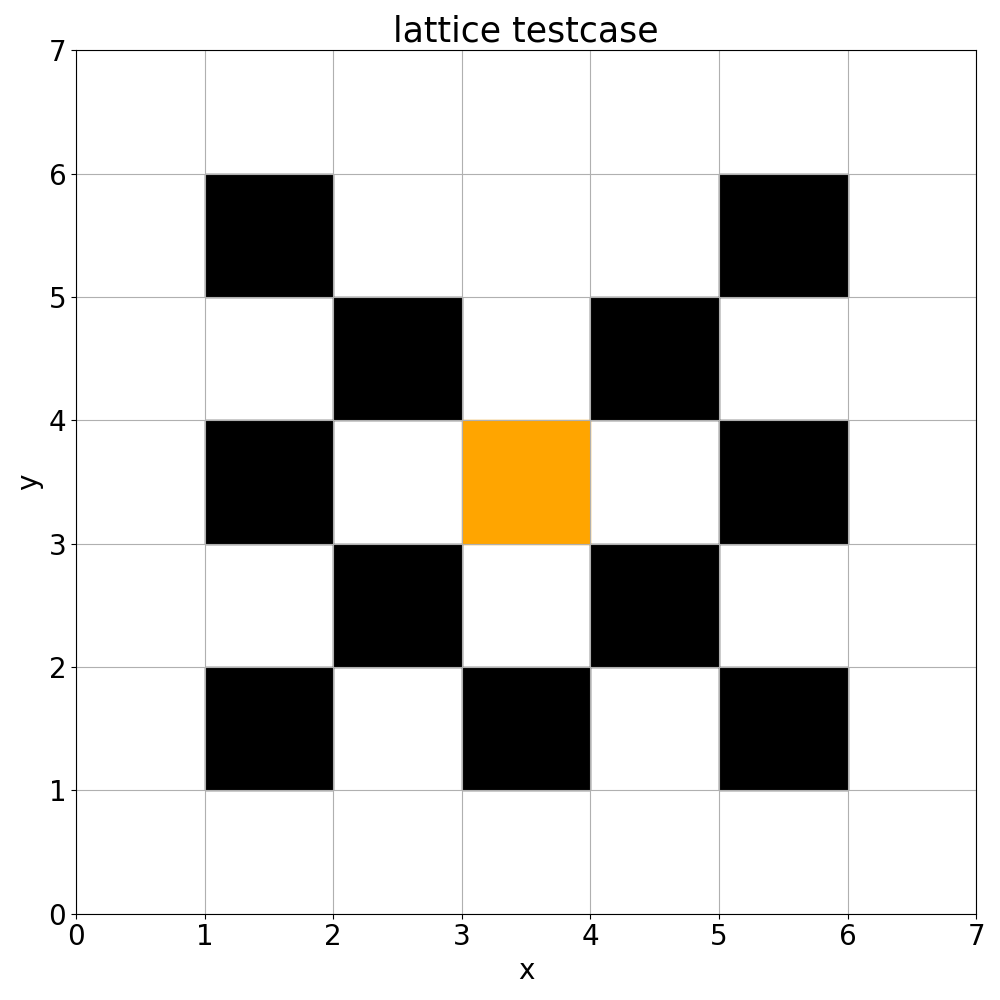}
         \caption{}\label{fig:setup_lattice_testcase}
     \end{subfigure}
     \hfill
     \begin{subfigure}[b]{0.49\textwidth}
         \centering
         \includegraphics[width=\textwidth]{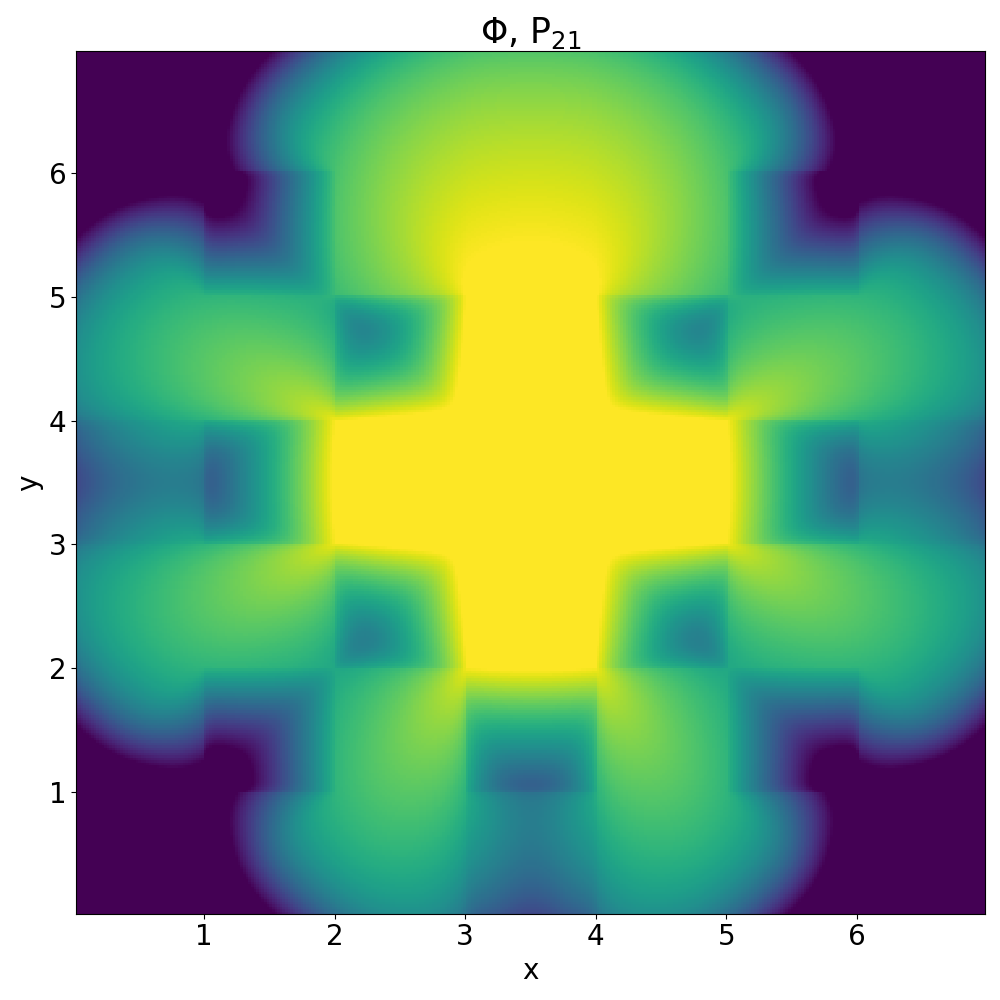}
         \caption{}
         \label{fig:scalar_flux_PN_Lattice_nx350_N21}
     \end{subfigure}
        \caption{
		Left: Setup lattice testcase. Right: P$_{21}$ solution of lattice testcase.}
	\label{fig:scalar_flux_PN_Lattice_nx350_N21_andsetup}
\end{figure}

\begin{figure}
	\includegraphics[width=0.49\textwidth]{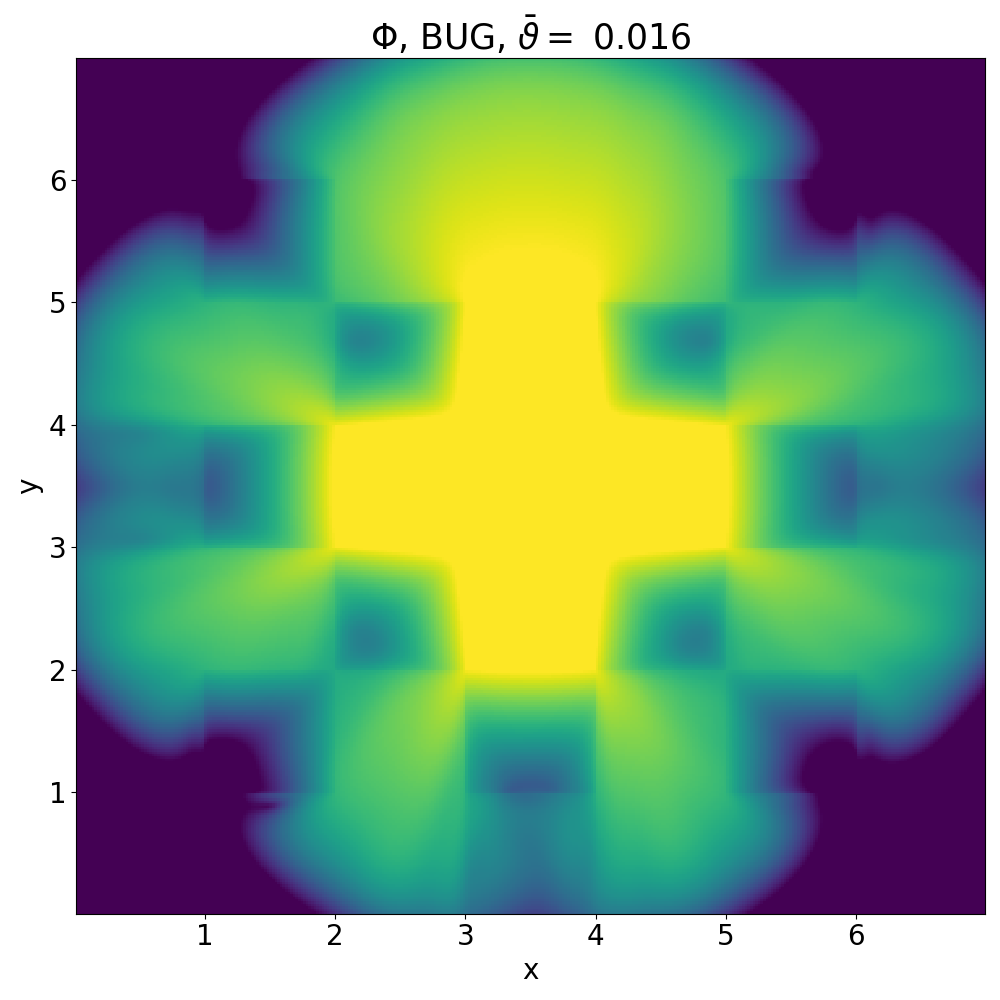}
	\includegraphics[width=0.49\textwidth]{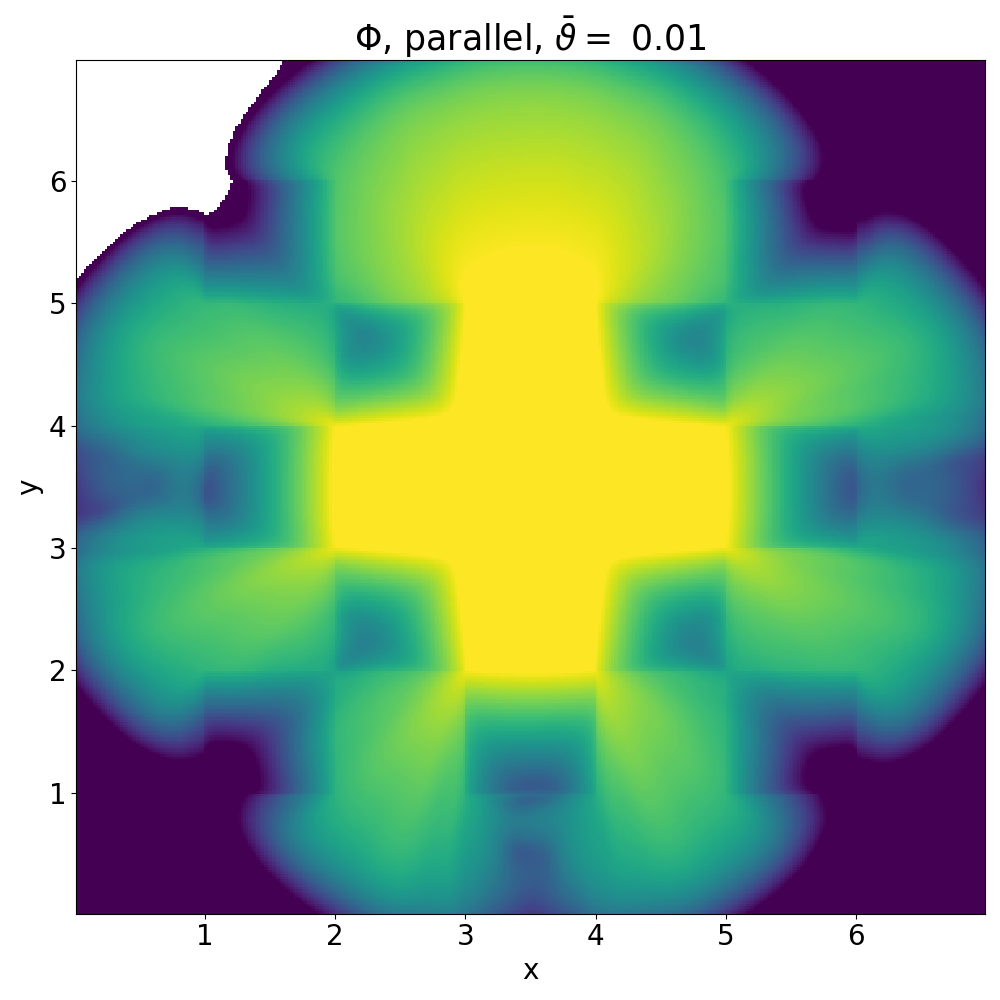}
    \includegraphics[width=0.49\textwidth]{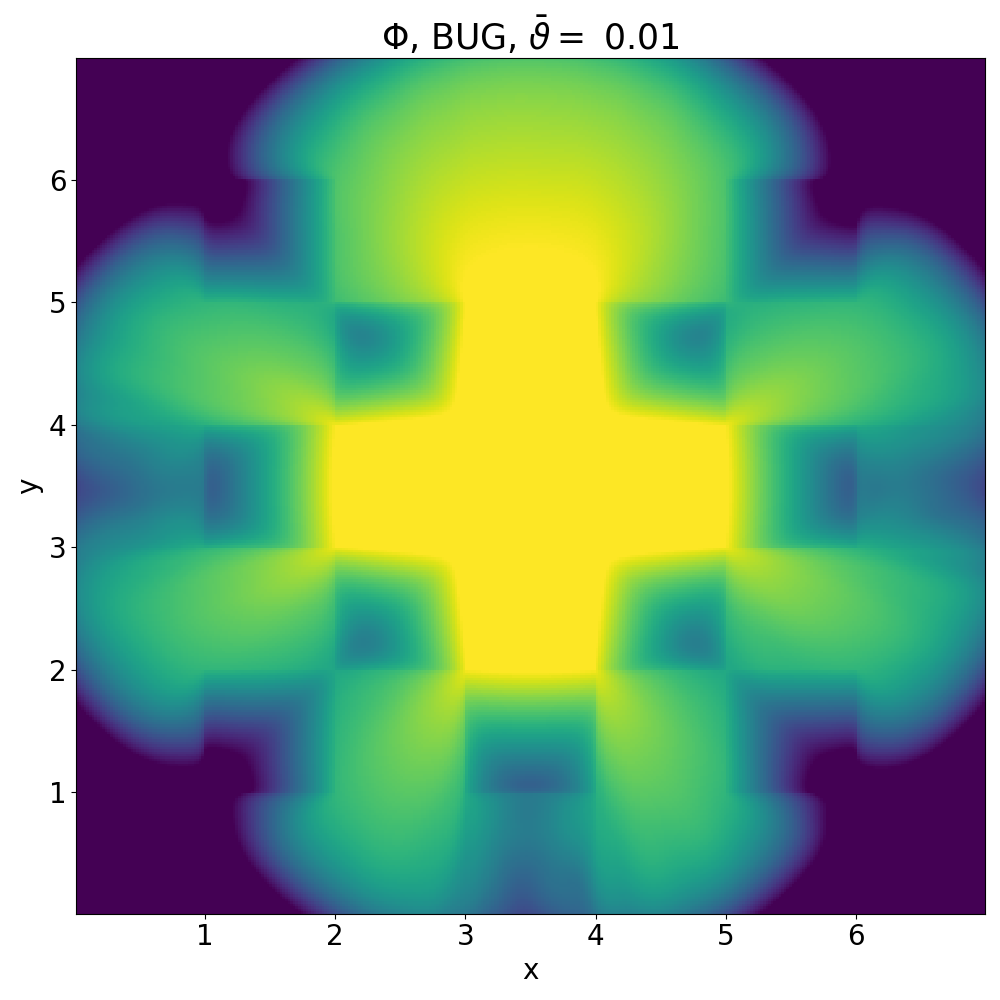}
	\includegraphics[width=0.49\textwidth]{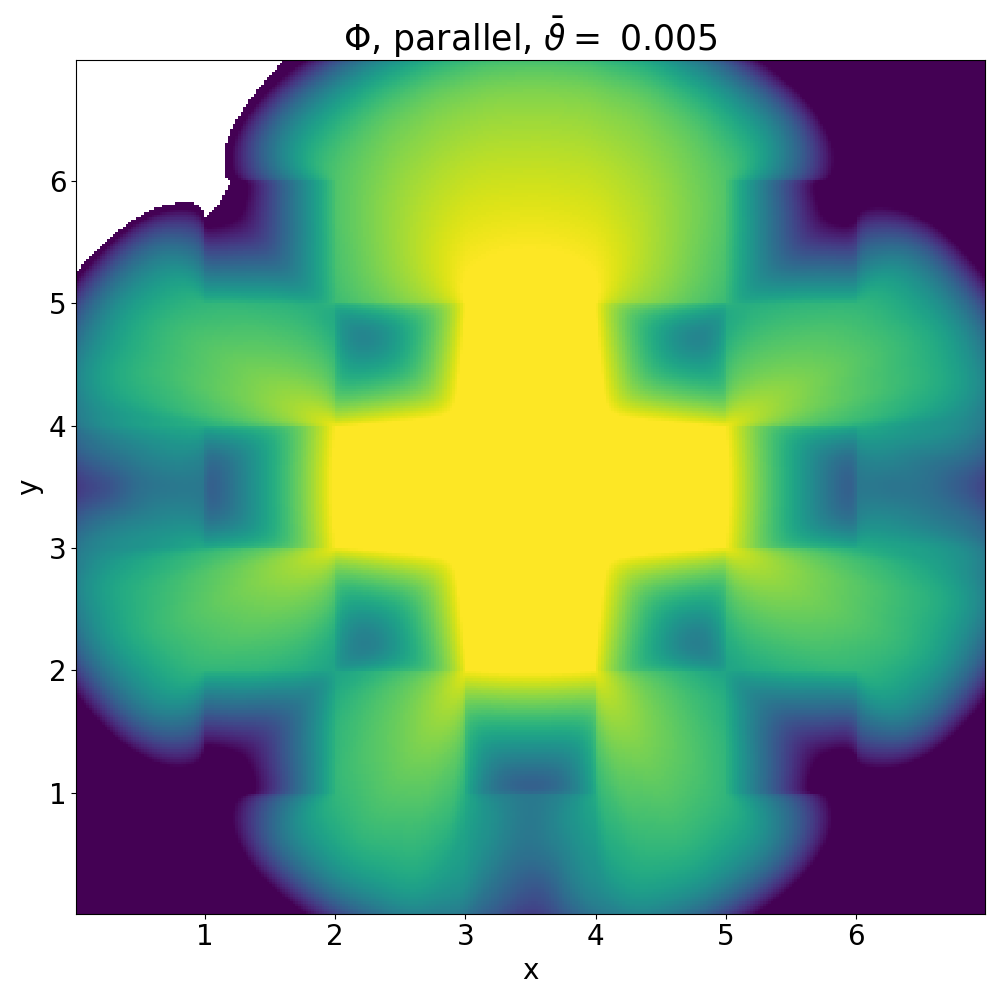}
	\caption{
		Left: Adaptive BUG integrator. Right: Parallel integrator. Negative values are depicted in white.}
	\label{fig:lattice_DLRA}
\end{figure}

In this section, we investigate the lattice testcase, which mimics the blocks in a nuclear reactor. For this, we wish to determine the solution to the two-dimensional radiative transfer equation
\begin{equation}
	\label{eq:rtela}
	\begin{aligned}
		&\partial_t f + \mathbf{\Omega}\cdot\nabla f + \sigma_t(\mathbf{x})f = \frac{\sigma_s(\mathbf{x})}{4\pi} \int_{-1}^{1} f \,d\mathbf{\Omega}+Q(\mathbf{x}),\qquad 
		(\mathbf{x}, \mathbf{\Omega}) \in [0,7]^2 \times \mathbb{S}^2 \;. \\[1mm]
		& f(t = 0) = 10^{-9}.
	\end{aligned}
\end{equation}
In contrast to previous numerical examples, the scattering and absorption cross-sections of the lattice testcase are spatially dependent. The spatial layout of the testcase can be found in Figure~\ref{fig:setup_lattice_testcase}. Here, white blocks indicate material which has only scattering properties, that is, $\sigma_s = 1$ and $\sigma_a = \sigma_t-\sigma_s = 0$. Black blocks are strongly absorbing with $\sigma_a = 10$ and no scattering, i.e., $\sigma_s = 0$. The orange block is an absorbing block with a source of $Q=1$. A reference solution has been computed with the spherical harmonics (P$_N$) method where we use a polynomial degree of $21$ which corresponds to $441$ expansion coefficients in the direction of flight $\mathbf{\Omega}$. The chosen number of spatial cells is $350$ in each spatial dimension, i.e., a total number of $350^2 = 122500$ cells is used. The scalar flux $\Phi(t,\mathbf x) = \int_{\mathbb{S}^2} f(t,\mathbf x,\mathbf{\Omega}) \,\mathrm{d}\mathbf{\Omega}$ is depicted in Figure~\ref{fig:scalar_flux_PN_Lattice_nx350_N21}. As commonly done for the lattice benchmark, we us a logarithmic scale.
In the following, we wish to approximate the P$_{21}$ solution with rank-adaptive dynamical low-rank approximation integrators. We recall that the rank adaptive BUG integrator will pick increased ranks compared to the parallel integrator when both methods use the same tolerance parameter $\vartheta$. Therefore, we run the BUG integrator with tolerances $\bar\vartheta=\{0.016, 0.01\}$. To obtain the same maximal rank, we run the parallel integrator with tolerances $\bar\vartheta=\{0.01, 0.005\}$. The corresponding scalar fluxes can be seen in Figure~\ref{fig:lattice_DLRA}.
\begin{figure}
\centering
	\includegraphics[width=0.75\textwidth]{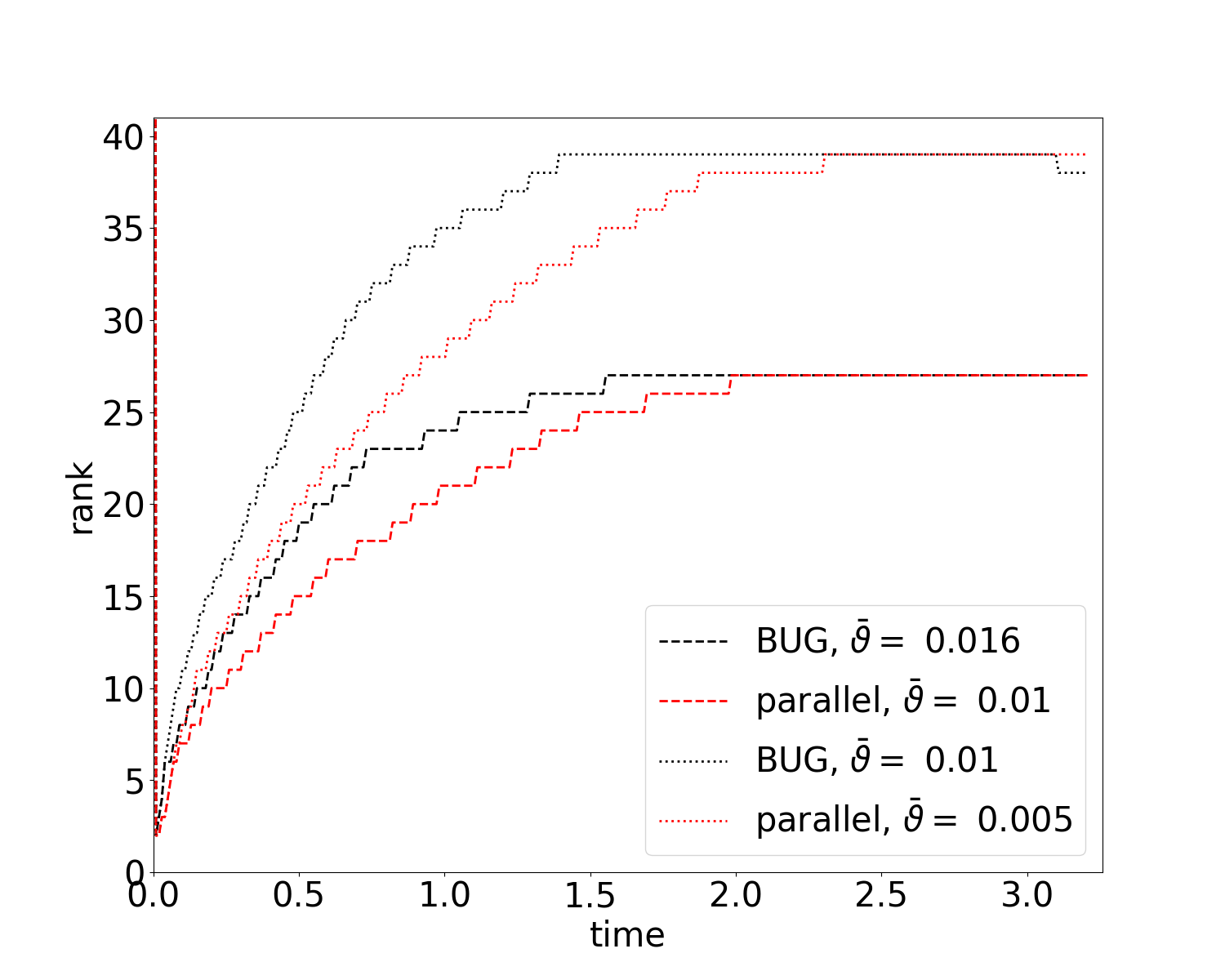}
	\caption{Rank evolution in time.}
	\label{fig:rank_lattice}
\end{figure}
It is observed that the parallel integrator and the BUG integrator yield similar results when the tolerance is chosen accordingly. However, the parallel integrator leads to negative values in certain regions which are depicted in white. This behaviour is not uncommon for modal approximations in radiative transfer and can also appear for the BUG integrator. The ranks chosen during the simulation, which can be found in Figure~\ref{fig:rank_lattice} are similar for both the BUG and the parallel integrator. Even though the parallel integrator increases the rank more slowly than BUG, the overall solution quality appears to not suffer from the reduced ranks except for regions with negative scalar fluxes. At the coarse tolerances ($\bar\vartheta = 0.016$ for BUG and $\bar\vartheta = 0.01$ for the parallel integrator) a maximal rank of $27$ is picked, which leads to minor numerical artifacts in the scattering regions. At the fine tolerances ($\bar\vartheta = 0.01$ for BUG and $\bar\vartheta = 0.005$ for the parallel integrator) a maximal rank of $39$ is chosen. Here, the approximations agree well with the full P$_{21}$ solution. In terms of computational cost, the runtime of the parallel integrator is 180 seconds for $\bar\vartheta = 0.01$ and 330 seconds for $\bar\vartheta = 0.005$. The BUG integrator has an increased runtime of 286 seconds for $\bar\vartheta = 0.016$ and 524 seconds for $\bar\vartheta = 0.01$. It is important to point out that the implementation of both method does not make use of the opportunity to compute $K$, $L$ (and for the parallel integrator $S$) in parallel, which would further favor the parallel integrator. The improved runtime of the parallel integrator results mostly from omitting the $2r \times 2r$ coefficient update.
\begin{figure}
\centering
	\includegraphics[width=0.49\textwidth]{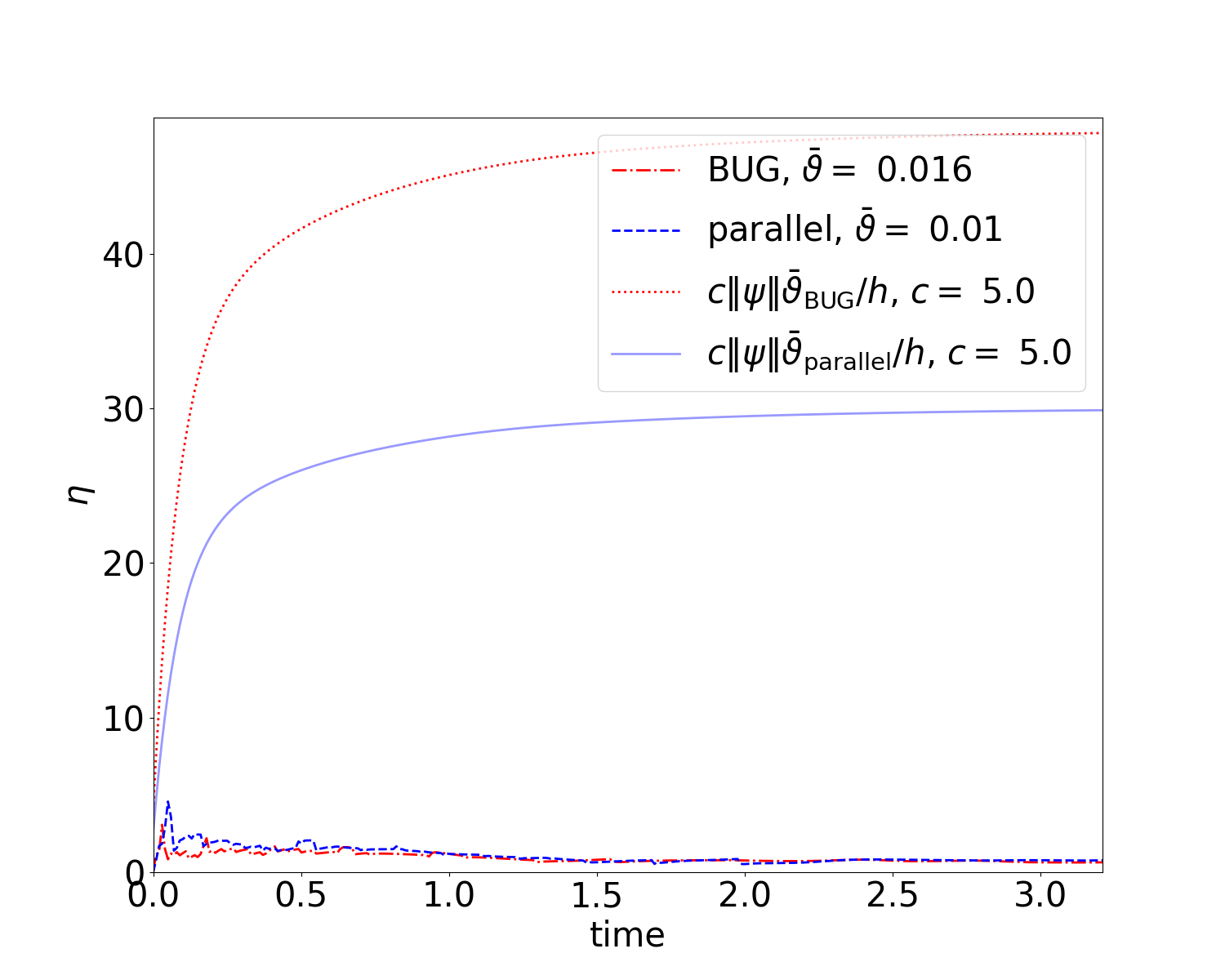}
 \includegraphics[width=0.49\textwidth]{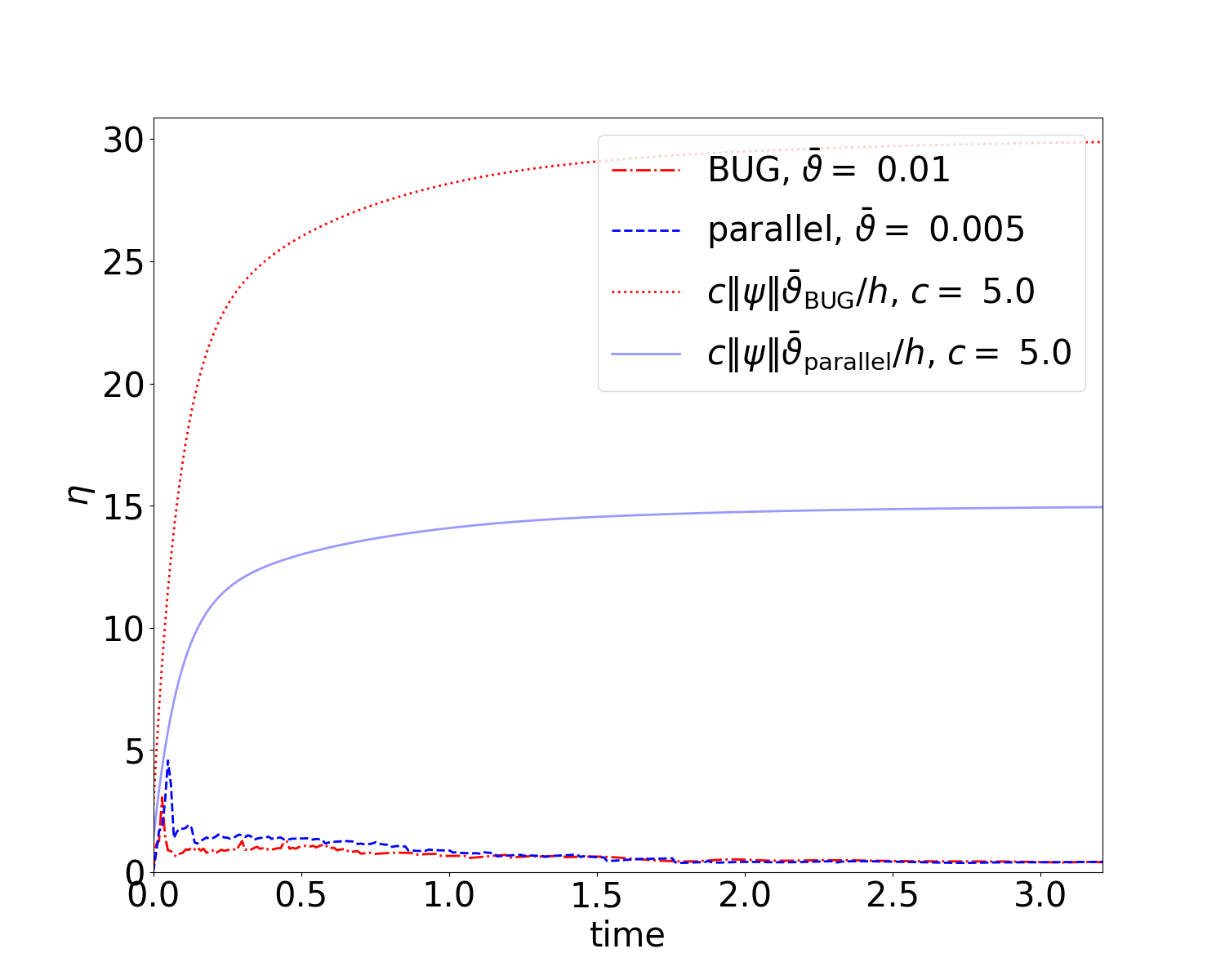}
	\caption{Values of $\eta$ in time.}
	\label{fig:eta_lattice}
\end{figure}

\subsection{Electron therapy}
\begin{figure}
     \centering
     \begin{subfigure}[b]{0.49\textwidth}
         \centering
         \includegraphics[width=\textwidth]{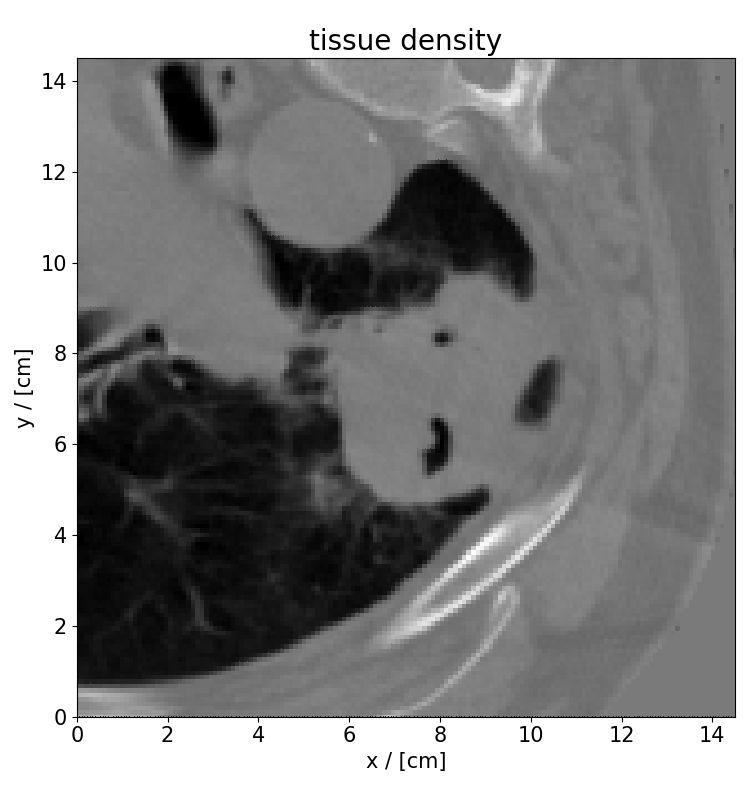}
         \caption{}\label{fig:setup_electron_therapy}
     \end{subfigure}
     \hfill
     \begin{subfigure}[b]{0.49\textwidth}
         \centering
         \includegraphics[width=\textwidth]{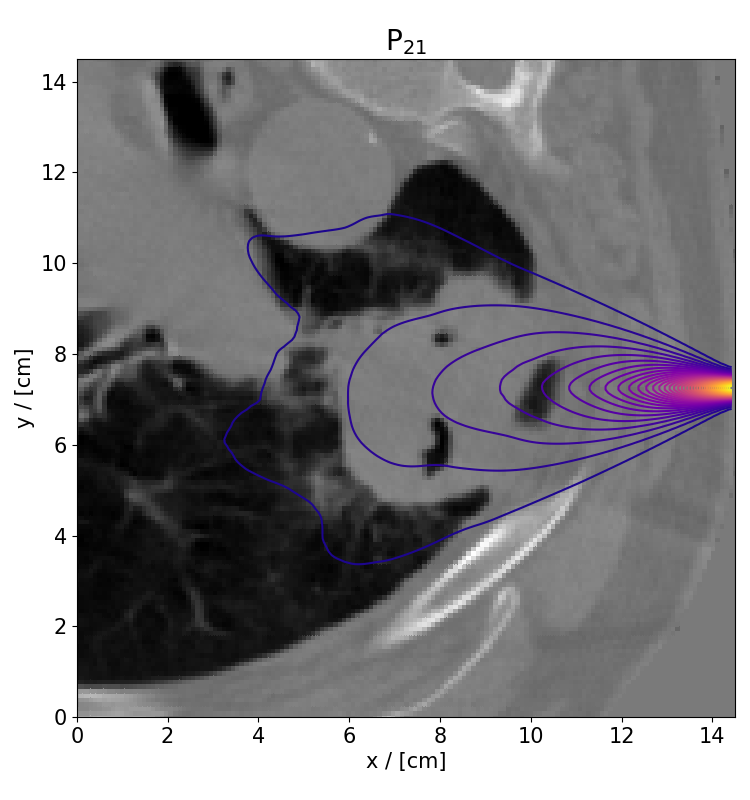}
         \caption{}
         \label{fig:reference_electron_therapy}
     \end{subfigure}
        \caption{
		Left: Setup radiation therapy application. Right: Corresponding P$_{21}$ solution.}
	\label{fig:dosesetupandfull}
\end{figure}

\begin{figure}[h!]
	\includegraphics[width=0.99\textwidth]{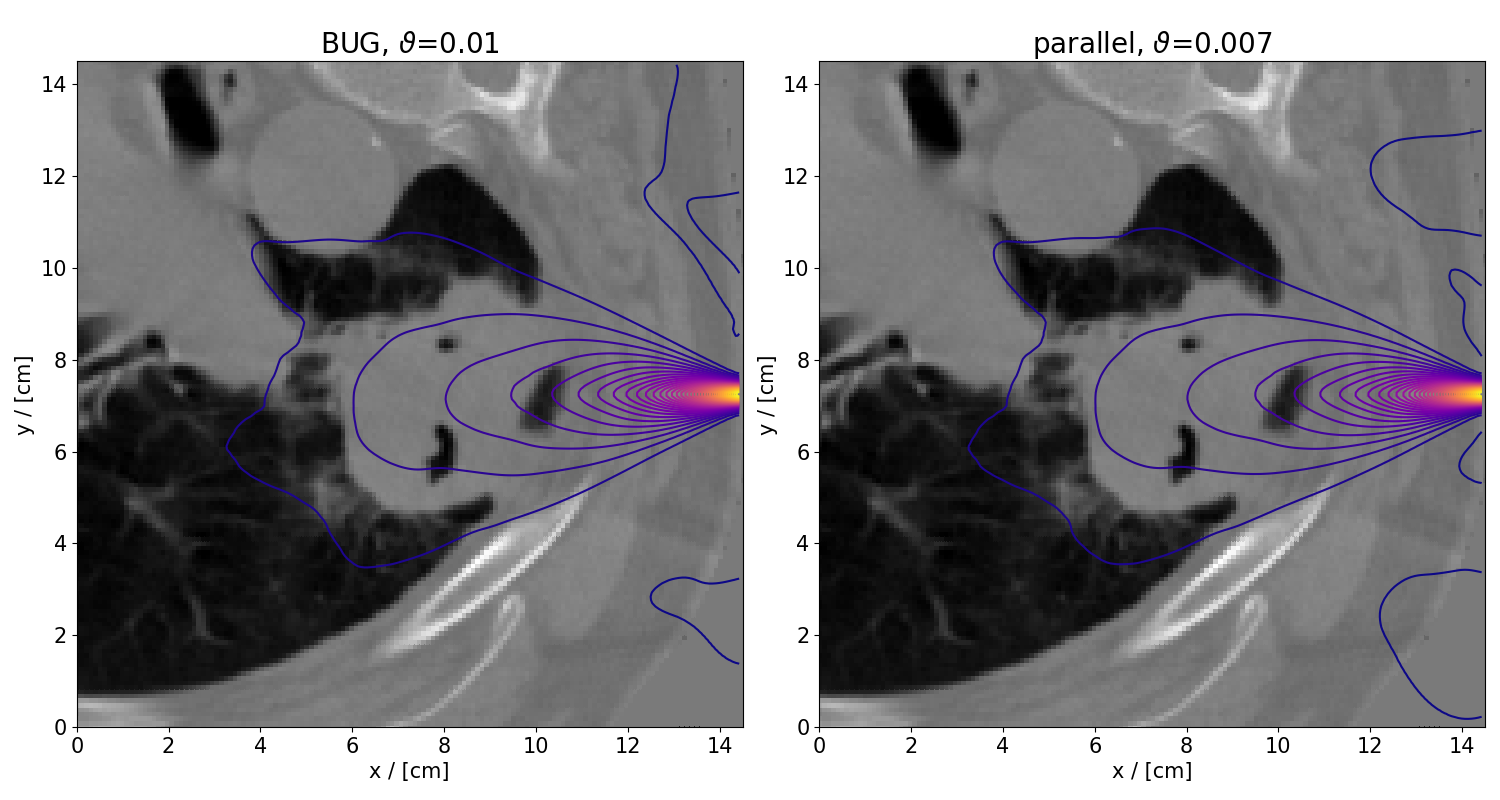}
	\caption{
		Left: Adaptive BUG integrator. Right: Parallel integrator.}
	\label{fig:rad_therapy_parallel_BUG}
\end{figure}

\begin{figure}
	\includegraphics[width=0.8\textwidth]{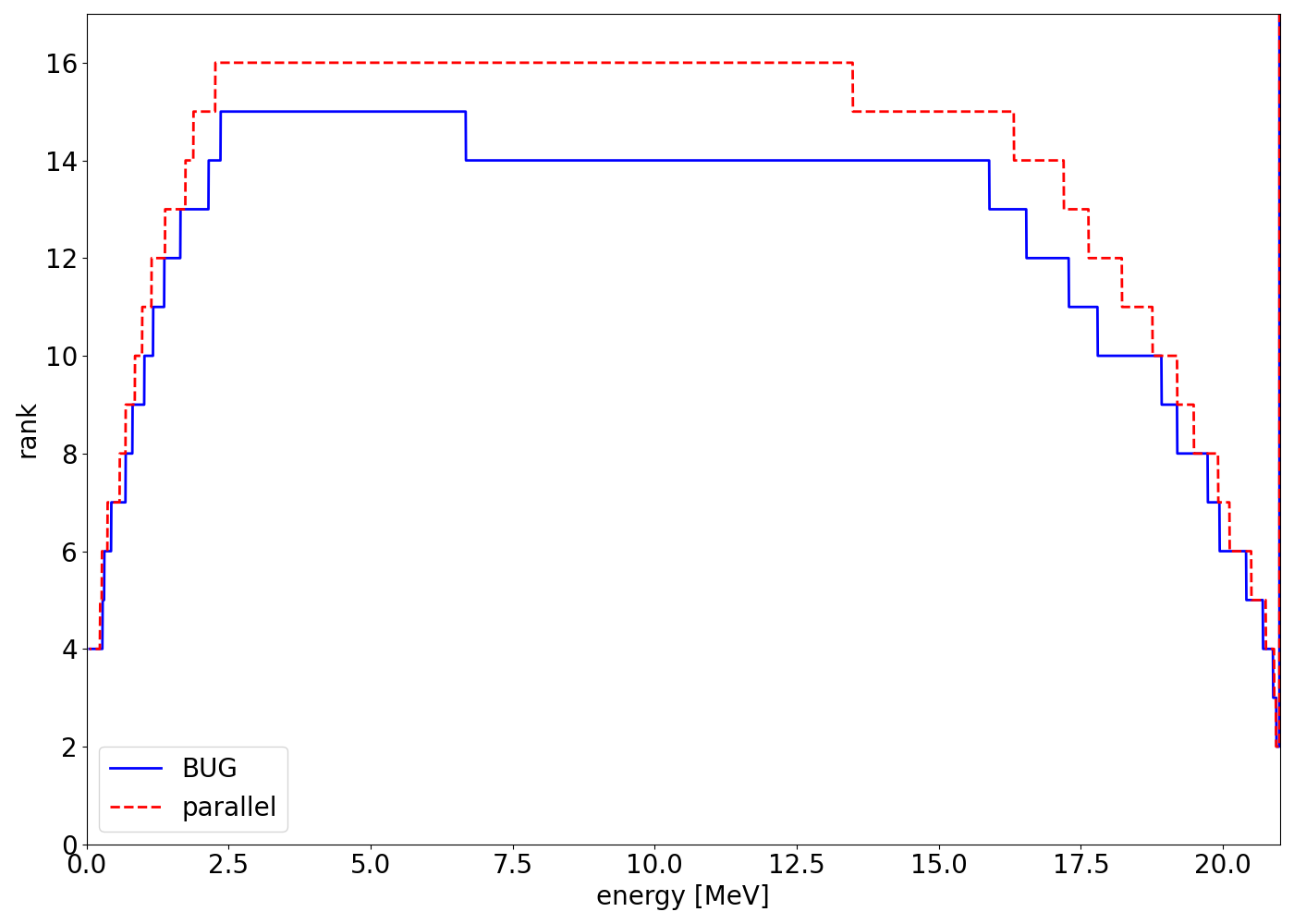}
	\caption{
		Rank over energy}
	\label{fig:ranksRadTherapy}
\end{figure}

\begin{figure}[h!]
	\includegraphics[width=0.99\textwidth]{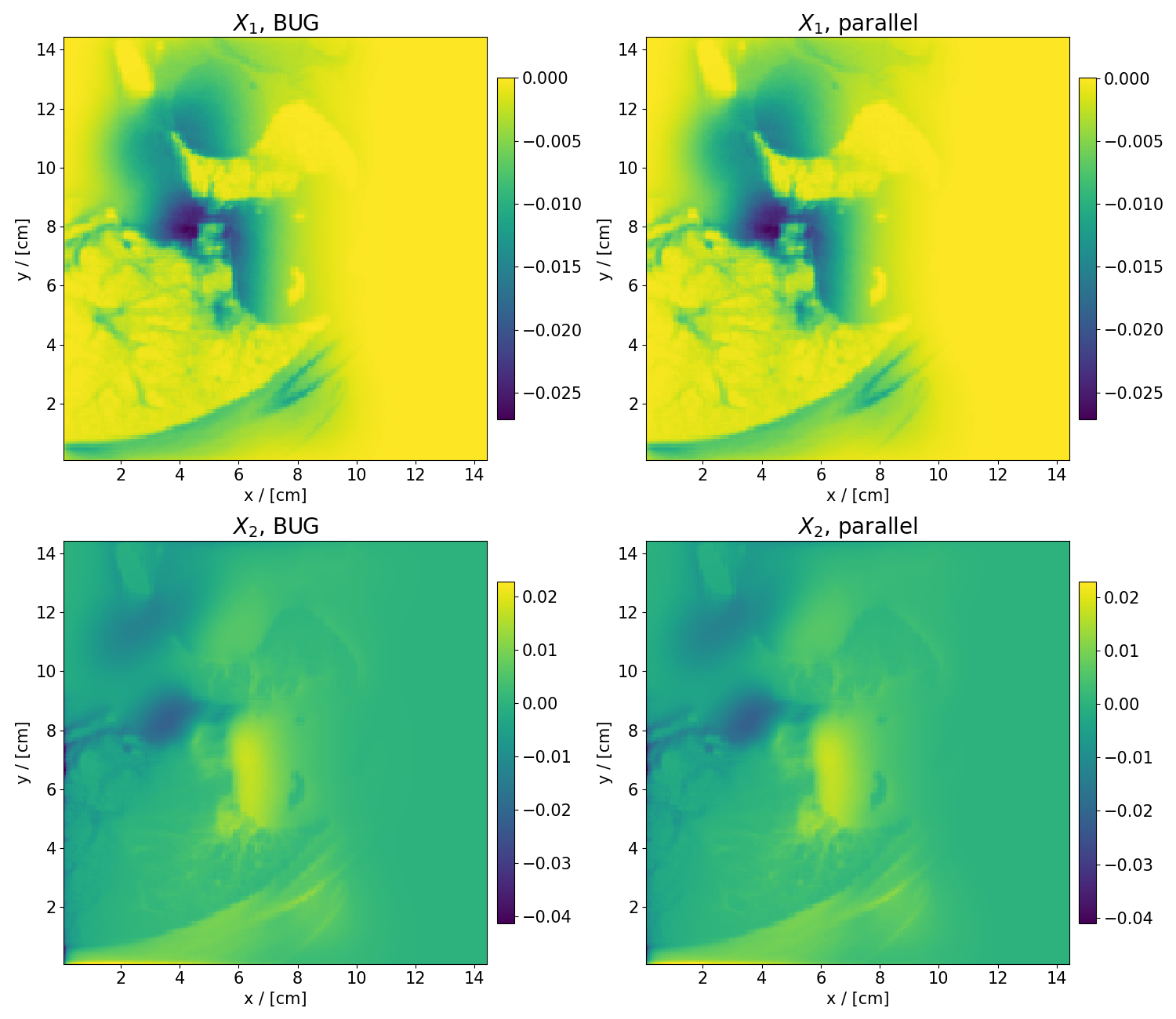}
	\caption{
		First two dominant spatial modes $X_1$ and $X_2$ at $E = 0$.}
	\label{fig:dominantXmodes}
\end{figure}

\begin{figure}[h!]
	\includegraphics[width=0.99\textwidth]{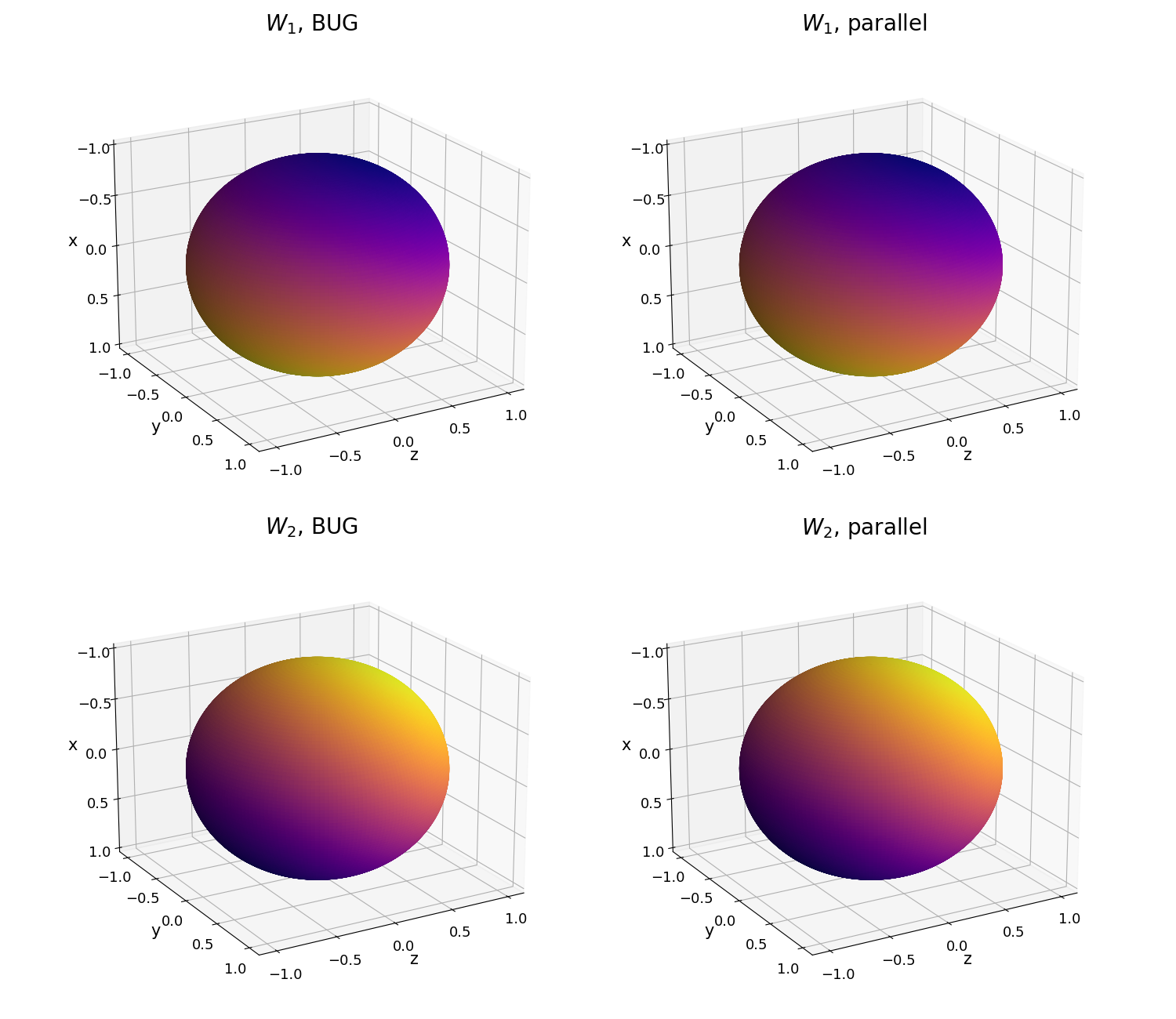}
	\caption{
		First two dominant directional modes $W_1$ and $W_2$ at $E = 0$.}
	\label{fig:dominantWmodes}
\end{figure}

In this section, we investigate the effectiveness of the parallel integrator for electron therapy. The main interest of computational radiation therapy is to predict the effectiveness of treatment plans by simulating the treatment procedure. The success of radiation treatment is determined by investigating the radiation dose absorbed by the patient at the tumor. At the same time, the radiation dose should not harm risk regions such as important organs or the spine. The movement of electrons through tissue is modeled by the continuous-slowing down equation 
\begin{align*}
    \mathbf{\Omega}\cdot\nabla\psi(E,\mathbf{x},\mathbf{\Omega})+\rho(\mathbf{x})\Sigma_t(E)\psi(E,\mathbf{x},\mathbf{\Omega}) = &\int_{\mathbb{S}^2}\rho(\mathbf{x})\Sigma_s(E,\mathbf{\Omega}\cdot\mathbf{\Omega}')\psi(E,\mathbf{x},\mathbf{\Omega}')\,d\mathbf{\Omega}'\nonumber\\
    &+\partial_E\left(\rho(\mathbf{x})S_t(E)\psi(E,\mathbf{x},\mathbf{\Omega})\right)\;,\\
     \psi(E,\mathbf{x},\mathbf{\Omega}) &= \psi_{\mathrm{BC}}(E,\mathbf{x},\mathbf{\Omega}) \quad \text{ for } \mathbf{x}\in\partial D\;,\\
     \psi(E_{\mathrm{max}},\mathbf{x},\mathbf{\Omega}) &= \psi_{\mathrm{max}}(\mathbf{x},\mathbf{\Omega})\;.
\end{align*}
In contrast to previous examples, we are interested in computing a steady state solution. Additionally, the radiation flux $\psi$ depends on the energy $E\in [0,E_{\mathrm{max}}]$. The tissue density $\rho(\mathbf{x})$ can be determined from a CT scan of the patient and the energy and angle-dependent cross sections $\Sigma_t$ and $\Sigma_s$ as well as the stopping power $S_t$ can be determined from physical data bases. To construct an efficient numerical algorithm, the energy is transformed into a pseudo-time $t$ according to
\begin{align}\label{eq:TildeE}
    t(E) := \int_0^{E_{\mathrm{max}}} \frac{1}{S_t(E')}\,dE'-\int_0^{E} \frac{1}{S_t(E')}\,dE'\;.
\end{align}
Then, for the transformed radiation flux
\begin{align}
    f(t,\mathbf{x},\mathbf{\Omega}):= S_t(E(t))\rho(\mathbf{x})\psi(E(t),\mathbf{x},\mathbf{\Omega})\;
\end{align}
we obtain the transformed continuous-slowing down equation  
\begin{align}\label{eq:CSD2}
\partial_t f = -\mathbf{\Omega}\cdot\nabla_{\mathbf{x}} \frac{f}{\rho}-\widetilde\Sigma_tf + \int_{\mathbb{S}^2}\widetilde\Sigma_s(t,\mathbf{\Omega}\cdot\mathbf{\Omega}')f(t,\mathbf{x},\mathbf{\Omega}')\,d\mathbf{\Omega}'\;,
\end{align}
where the tilde denotes a transformation of energy. We then evolve the solution in (pseudo) time $t$ with the BUG and parallel integrators. For more details on the chosen DLRA approximation we refer to \cite{kusch2021robust}. To test the parallel integrator, we follow \cite{kusch2021robust} and determine the dose 
\begin{align}
    D(\mathbf{x})=\frac{1}{\rho(\mathbf{x})}\int_0^{\infty}\int_{\mathbb{S}^2}S_t(E)\psi(E,\mathbf{x},\mathbf{\Omega})\,d\mathbf{\Omega} dE\;
\end{align}
when radiating a lung tumor with an energy beam of strength $E_{\mathrm{max}} = 21$ which is modelled by
\begin{linenomath*}\begin{align*}
    \psi_{\text{in}}(E,\mathbf x, \mathbf \Omega) = 10^5&\cdot\exp(-(\Omega_{1,*}-\Omega_1)^2/\sigma_{\Omega_1})\cdot\exp(-(E_{\text{max}}-E)^2/ \sigma_E)\\
    &\cdot\exp(-(x_*-x)^2/\sigma_x)\cdot\exp(-(y_*-y)^2/\sigma_y)\;.
\end{align*}\end{linenomath*}
Parameters of the ingoing beam are \\[1mm]
\begin{center}
    \begin{tabular}{ | l | p{8cm} |}
    \hline
    $x_* = 7.25, y_* = 14.5$ & spatial mean of particle beam in cm\\
    $\Omega_{1,*} = 1$ & directional mean of particle beam \\
    $\sigma_{\Omega_1}^{-1}= 75$ & inverse directional beam variance \\
    $\sigma_x^{-1}=\sigma_y^{-1} = 20$ & inverse spatial beam variance \\
    $\sigma_E^{-1}=100$ & inverse energy variance \\
    \hline
    \end{tabular}
\end{center}
We use the upwind method from previous examples with an equidistand spatial grid of $200 \times 200$ spatial cells. Again a modal discretization of angle is chosen, where we use spherical harmonics basis functions of degree up to $21$. The tissue density row is depicted in Figure~\ref{fig:setup_electron_therapy} and the corresponding full rank P$_{21}$ solution which is depicted in Figure~\ref{fig:reference_electron_therapy} has been computed with a runtime of 17467 seconds.
To reduce computational costs we approximate the dose distribution with DLRA using the BUG integrator with a tolerance of $\vartheta = 0.01$ and compare it to the parallel integrator with a tolerance of $\vartheta = 0.007$, which yields the comparable maximal rank. The corresponding doses are shown in Figure~\ref{fig:rad_therapy_parallel_BUG}. We observe that both integrators yield satisfactory approximations which agree well with the P$_{21}$ reference solution. The runtimes are 5423 seconds for the BUG integrator and 3894 seconds for the parallel integrator. As expected, we observe fluctuations on the upper and lower left side, which are negligibly small. The corresponding ranks over energy can be found in Figure~\ref{fig:ranksRadTherapy}. Dominant spatial modes computed with BUG and the parallel integrator can be found in Figure~\ref{fig:dominantXmodes} and dominant directional modes are depicted in Figure~\ref{fig:dominantWmodes}.

